\RequirePackage{pdf14}

\RequirePackage{fix-cm}
\documentclass[smallextended]{svjour3}       % onecolumn (second format)
\smartqed  % flush right qed marks, e.g. at end of proof
\usepackage{graphicx}

\usepackage[T1]{fontenc}

\usepackage{amssymb,latexsym,amsmath}
\usepackage{mathrsfs}
\usepackage{enumitem}
\usepackage{caption}
\usepackage{comment}
\usepackage{inputenc}
 \usepackage{makeidx}
\usepackage{color}

\allowdisplaybreaks

\allowdisplaybreaks

\usepackage{epsfig} % for postscript graphics files
\usepackage{verbatim}

\def\qed{\hfill $\diamond$}

\usepackage{amsfonts}
\usepackage{float}
\usepackage{multirow}
       \newtheorem{assumption}{\bf{Assumption}}%[section]
\usepackage{color}%

\def\sY{{\mathbb Y}}
\allowdisplaybreaks

\renewcommand{\qed}{\hfill$\square$}

\newcommand{\sy}[1]{{\color{black} #1}}

\begin{document}

%\begin{frontmatter}
%\allowdisplaybreaks

\title{Another Look at Partially Observed Optimal Stochastic Control: Existence, Ergodicity, and Approximations without Belief-Reduction}

\titlerunning{Partially Observed Optimal Stochastic Control without Belief-Reduction}        % if too long for running head

%\title{Partially Observed MDPs without Belief Reduction: Existence, Invariance, and Approximation Results}
%\title{Another Look at Optimal Control under Partial Information without Belief-MDP Reduction and conditionally-exogenous variable-independent Policies}

\author{Serdar Y\"uksel}

\institute{S. Y\"{u}ksel \at
              Department of Mathematics and Statistics \\
              Queen's University \\
              \email{yuksel@queensu.ca}. This research was partially supported by the Natural Sciences and Engineering Research Council of Canada (NSERC).           %  \\
%             \emph{Present address:} of F. Author  %  if needed
}

%\thanks{Department of Mathematics and
%    Statistics, Queen's University, Kingston, Ontario, Canada, K7L
%    3N6.  Email: yuksel@mast.queensu.ca. This research was
%    partially supported by the Natural Sciences and Engineering
%    Research Council of Canada (NSERC).}
%}
%%\begin{document}

\maketitle

\begin{abstract}
We present an alternative view for the study of optimal control of partially observed Markov Decision Processes (POMDPs). We first revisit the traditional (and by now standard) separated-design method of reducing the problem to fully observed MDPs (belief-MDPs), and present conditions for the existence of optimal policies. Then, rather than working with this standard method, \sy{we define a Markov chain taking values in an infinite dimensional product space with the history process serving as the controlled state process and a further refinement in which the control actions and the state process are causally conditionally independent given the measurement/information process}. 
We provide new sufficient conditions for the existence of optimal control policies under the discounted cost and average cost infinite horizon criteria. In particular, while in the belief-MDP reduction of POMDPs, weak Feller condition requirement imposes total variation continuity on either the system kernel or the measurement kernel, with the approach of this paper only weak continuity of both the transition kernel and the measurement kernel is needed (and total variation continuity is not) together with regularity conditions related to filter stability. For the discounted cost setup, we establish near optimality of finite window policies via a direct argument involving near optimality of quantized approximations for MDPs under weak Feller continuity, where finite truncations of memory can be viewed as quantizations of infinite memory with a uniform diameter in each finite window restriction under the product metric. For the average cost setup, we provide new existence conditions and also a general approach on how to initialize the randomness which we show to establish convergence to optimal cost.  In the control-free case, our analysis leads to new and weak conditions for the existence and uniqueness of invariant probability measures for nonlinear filter processes, where we show that unique ergodicity of the measurement process and a measurability condition related to filter stability leads to unique ergodicity.
\end{abstract}

%\begin{keyword}
%{\bf Key words}: Partially Observed MDPs, nonlinear filtering
%%\kwd{\LaTeXe}
%\end{keyword}

%\begin{keyword} 
%\end{keyword}

%\end{frontmatter}

%{\cal P}agestyle{myheadings}
%\thispagestyle{plain}

\section{Introduction, Literature Review, and Main Results}

Partially observed Markov Decision processes (POMDPs) present challenging mathematical problems with significant applied relevance.  It is known that any POMDP can be reduced to a (completely observable) MDP \cite{Yus76}, \cite{Rhe74}, whose states are the posterior state distributions or {\it beliefs} of the observer; in particular, the belief-MDP is a (fully observed) Markov decision process. 

%Accordingly, for finite horizon problems and a large class of infinite horizon discounted cost problems it is a standard result that an optimal control policy will use the belief as a sufficient statistic for optimal policies (see \cite{Yus76,Rhe74,Blackwell2}). Hence, the POMDP and the corresponding belief-MDP are equivalent in the sense of cost minimization. %Therefore, results developed for MDPs can be applied to the belief-MDP.

However, the use of belief-MDPs requires one to establish several regularity, continuity and stability (such as filter stability or unique ergodicity) results for arriving at existence, finite memory approximations, robustness, as well as learning theoretic results (see e.g. \cite{yu2008near,kara2020near,kara2021convergence,MYRobustControlledFS}. While significant progress on the regularity properties of belief-MDPs has been made in the literature, in this paper we will see that further refinements are possible if one does not restrict the analysis to a belief-MDP based formulation. In this paper, we present an alternative view for the study of infinite horizon average-cost or discounted cost optimal control of POMDPs. Our approach will be on a design which is not based on separation / or belief-MDP reduction. We define a Markov chain taking values in an infinite dimensional product space with control actions and the state process causally conditionally independent given the measurement process. 

For the controlled case, we provide new sufficient conditions for the existence of optimal control policies for average and discounted cost criteria, and under the latter criterion our analysis will establish a general result on near optimality of finite memory policies. In the control-free case, our analysis leads to sufficient conditions for the existence and uniqueness of an invariant probability measure for nonlinear filters. 

We now present the problem. Consider a stochastic process $\{X_t, t \in \mathbb{Z}_+\}$, where each element $X_t$ takes values in some standard Borel space $\mathbb{X}$, with dynamics described by
\begin{eqnarray}\label{updateEq}
 X_{t+1} &=& F(X_t, U_t, W_t) \label{updateEq1} \\
Y_{t} &=& G(X_t, V_t) \label{updateEq2} 
\end{eqnarray}
where $Y_t$ is an $\mathbb{Y}$-valued measurement sequence; we take $\mathbb{Y}$ also to be some standard Borel space. Suppose further that $X_0 \sim \mu$. Here, $W_t, V_{t}$ are mutually independent i.i.d. noise processes. This system is subjected to a control/decision process where the control/decision at time $n$, $U_n$, incurs a cost $c(X_n,U_n)$. The decision maker only has access to the measurement process $Y_n$ and $U_n$ causally: An {\em admissible policy} $\gamma$ is a
sequence of control/decision functions $\{\gamma_t,\, t\in \mathbb{Z}_+\}$ such
that $\gamma_t$ is measurable with respect to the $\sigma$-algebra
generated by the information variables
\[
I_t=\{Y_{[0,t]},U_{[0,t-1]}\}, \quad t \in \mathbb{N}, \quad
  \quad I_0=\{Y_0\}.
\]
so that
\begin{equation}
\label{eq_control}
U_t=\gamma_t(I_t),\quad t\in \mathbb{Z}_+
\end{equation}
are the $\mathbb{U}$-valued, where $\mathbb{U}$ is assumed to be a compact action space, control/decision
actions and we use the notation
\[
Y_{[0,t]} = \{Y_s,\, 0 \leq s \leq t \}, \quad U_{[0,t-1]} =
  \{U_s, \, 0 \leq s \leq t-1 \}.
\]
We define $\Gamma$ to be the set of all such (strong-sense) admissible policies. We emphasize the implicit assumption here that the control policy can also depend on the prior probability measure $\mu$. 

\sy{We assume that all of the random variables are defined on a common probability space $(\Omega, {\cal F}, P)$}. We note that (\ref{updateEq1})-(\ref{updateEq2}) can also, equivalently (via stochastic realization results \cite[Lemma~1.2]{gihman2012controlled} \cite[Lemma~3.1]{BorkarRealization}, \cite[Lemma F]{aumann1961mixed}), be represented with transition kernels: the state transition kernel is denoted with ${\cal T}$ so that for Borel $B \subset \mathbb{X}$ \[{\cal T}(B|x,u) := P(X_1 \in B | X_0=x,U_0=u). \quad\] We will denote the measurement kernel with $Q$ so that for Borel $B \subset \mathbb{Y}$: \[Q(B|x) := P(Y_0 \in B | X_0=x).\]

\sy{Let $c: \mathbb{X} \times \mathbb{U} \to \mathbb{R}_+$ be a stage-wise cost function. Throughout we will assume that $c$ is continuous and bounded.} For (\ref{updateEq})-(\ref{updateEq2}), we are interested in minimizing either the average-cost optimization criterion
\begin{eqnarray}\label{expCost}
J_{\infty}(\mu,\gamma) := \limsup_{N \to \infty} {1 \over N} E^{\gamma}_{\mu}[\sum_{t=0}^{N-1} c(X_t, U_t)] 
\end{eqnarray}
or the discounted cost criterion (for some $\beta \in (0,1)$
\begin{eqnarray}\label{expDiscCost}
J_{\beta}(\mu,\gamma) :=E^{\gamma}_{\mu}[\sum_{t=0}^{\infty} \beta^k c(X_t, U_t)] 
\end{eqnarray}
%or sample path costs
%\begin{eqnarray}\label{sampCost}
%\limsup_{N \to \infty} {1 \over N} \sum_{t=0}^{N-1} c(X_t, U_t)
%\end{eqnarray}
over all admissible control policies $\gamma = \{\gamma_0, \gamma_1, \cdots,\} \in \Gamma$ with $X_0 \sim \mu$. \sy{Here $E^{\gamma}_{\mu}[\cdot]$ refers to the expectation given policy $\gamma$ and initial measure $\mu$}. With ${\cal P}(\mathbb{U})$ denoting the set of probability measures on $\mathbb{U}$ endowed with the weak convergence topology, we will also, when needed, allow for independent randomizations so that $\gamma_n(I_n)$ is ${\cal P}(\mathbb{U})$-valued for each realization of $I_n$.  

For such average cost problems, in some applications, an optimality result may only hold for a restrictive class of initial conditions or initializations. This will be noted more explicitly later in the paper.

We will also consider the control-free case where the system equation (\ref{updateEq1}) does not have control dependence; in this case only a decision is to be made at every time stage and $U$ is present only in the cost expression in (\ref{expCost}). This important special case has been studied extensively in the theory of nonlinear filtering.

\subsection{Literature review and preliminaries}
In the following, we present a brief literature review on optimal control of POMDPs, before presenting the main results of the paper. 

\noindent{\bf POMDPs, separated policies and belief-MDPs.} It is well-known that any POMDP can be reduced to a (completely observable) MDP \cite{Yus76}, \cite{Rhe74}, whose states are the posterior state probabilities, or beliefs, of the observer; that is, the state at time $t$ is
\begin{align}
\pi_t(\,\cdot\,) := P\{X_{t} \in \,\cdot\, | Y_0,\ldots,Y_t, U_0, \ldots, U_{t-1}\} \in {\cal P}(\mathbb{X}). \nonumber
\end{align}
We call this equivalent MDP the belief-MDP\index{Belief-MDP}. The belief-MDP has state space ${\cal P}(\mathbb{X})$ and action space $\mathbb{U}$. Here, ${\cal P}(\mathbb{X})$ is equipped with the Borel $\sigma$-algebra generated by the topology of weak convergence \cite{Bil99}. Since $\mathbb{X}$ is a Borel space, ${\cal P}(\mathbb{X})$ is metrizable with the Prokhorov metric which makes ${\cal P}(\mathbb{X})$ into a Borel space \cite{Par67}. The transition probability $\eta$ of the belief-MDP can be constructed as follows (see also \cite{Her89}). If we define the measurable function 
\[F(\pi,u,y) := P\{X_{t+1} \in \,\cdot\, | \pi_t = \pi, U_t = u, Y_{t+1} = y\}\]
 from ${\cal P}(\mathbb{X})\times\mathbb{U}\times\sY$ to ${\cal P}(\mathbb{X})$ and the stochastic kernel $H(\,\cdot\, | \pi,u) := P\{Y_{t+1} \in \,\cdot\, | \pi_t = \pi, U_t = u\}$ on $\sY$ given ${\cal P}(\mathbb{X}) \times \mathbb{U}$, then $\eta$ can be written as
\begin{align}
\eta(\,\cdot\,|\pi,u) = \int_{\sY} 1_{\{F(\pi,u,y) \in \,\cdot\,\}} H(dy|\pi,u). \label{kernelFilter}
\end{align}
The one-stage cost function $c$ of the belief-MDP is given by
\begin{align}
\tilde{c}(\pi,u) := \int_{\mathbb{X}} c(x,u) \pi(dx). \label{weak:eq8}
\end{align}
In particular, the belief-MDP is a (fully observed) Markov decision process with the components $({\cal P}(\mathbb{X}),\mathbb{U},\eta,\tilde{c})$. %For the belief-MDP define the history spaces $\tilde{\sH}_{t}=({\cal P}(\mathbb{X}) \times\mathbb{U})^{t}\times{\cal P}(\mathbb{X})$, $t=0,1,2,\ldots$ and let $\tilde{\Pi}$ denote the set of all policies for the belief-MDP, where the policies are defined in an usual manner. Let $\tilde{W}(\tilde{\Pi},\xi)$ denote the cost function of policy $\tilde{\Pi} \in \tilde{\Pi}$ for initial distribution $\xi$ of the belief-MDP, where $\tilde{W} \in \{J,V\}$.

For finite horizon problems and a large class of infinite horizon discounted cost problems it is a standard result that an optimal control policy will use the belief $\pi_t$ as a sufficient statistic for optimal policies (see \cite{Yus76,Rhe74,Blackwell2}).% This separated approach of reducing a POMDP to a belief-based MDP is remarkably powerful, in view of existence results, structural results, and also on approximation results \cite{SYLTAC2017POMDP} \cite{zhou2010solving}. Theorem \ref{weakFellerBeliefMDP}, to be reviewed below, presents conditions for which such a separated design can be facilitated via measurable selection conditions for fully observed MDPs, where the fully observed state is probability measure valued.

%However, non-separated design can also be used to arrive at existence results.  \cite[Section 5.4.2, Theorem 5.6]{YukselWitsenStandardArXiv} presents such a characterization, without reducing the problem to a belief-MDP, though the conditions obtained correspond to those given in Assumption \ref{TV_channel} since a reduction is to be imposed, requiring total variation continuity, as a result of the pointwise convergence of the density functions.

%It would be interesting to establish such ergodicity results for partially observed controlled models under a fixed stationary policy (stationary in the filter realization); while we will discuss this problem with an alternative approach for the control-free case, the paper does not study this problem in the controlled setup. 

Several papers \cite{platzman1980optimal,fernandez1990remarks,runggaldier1994approximations,hsu2006existence} study the average-cost control problem under the assumption that the state space is finite; they provide reachability type conditions for the belief kernels. Reference \cite{Bor00} considers the finite model setup and \cite{borkar2004further} considers the case with finite-dimensional real-valued state spaces under several technical conditions on the controlled state process and \cite{StettnerSICON19} considers several conditions directly on the filter process leading to an equi-continuity condition on the relative discounted value functions. In a related discussion, \cite{feinberg2012average} studies the existence problem for the average-cost control problem under weak continuity conditions for the controlled kernels. In all of the aforementioned studies, the vanishing discount method is considered. Further discussion will be presented later in the paper.

%This method almost always requires strong ergodicity conditions or recurrence conditions to allow for convergence of the discounted value functions to a limit function (e.g., to utilize the Arzela-Ascoli theorem, which would require equi-continuity type conditions, as the discount parameter $\beta$ approaches unity). %The ergodicity properties of the non-filter process itself is a challenging problem with limited results even in the control-free setup (as reviewed above), therefore it is remarkable that e.g. \cite{Bor00} has established existence results through the vanishing discount method. Nonetheless, we will observe that the conditions required can be further relaxed. Among efforts for obtaining approximate solutions for POMDPs, \cite{white1994finite} \cite{kara2020near} provide bounds and convergence results for the discounted infinite horizon case, where the controllers use finite memory. See \cite{SYLTAC2017POMDP} \cite{zhou2010solving} for a detailed literature analysis and further results.

%{\bf Regularity of Belief-MDPs.} 

\noindent{\bf Control-free setup.} For the special case without control, the belief process is known as the (nonlinear) filter process, and by the discussion above, this itself is a Markov process. The stability properties of such processes has been studied, where the existence of an invariant probability measure for the belief process, as well as the uniqueness of such a measure (i.e., the unique ergodicity property) has been investigated under various conditions, see. e.g. \cite{budhiraja2002invariant}. For the control-free case, \cite{chigansky2010complete} provides a comprehensive discussion on both the ergodicity of the filter process as well as filter stability, when the state space is finite but the measurement space is not necessarily so, under the further assumption that the unobserved state process is stationary and ergodic. \cite[Theorem 2]{DiMasiStettner2005ergodicity} and \cite[Prop 2.1]{van2009uniformSPA} assume that the hidden state process is ergodic and the filter is stable (almost surely or in expectation under total variation); these papers crucially embed the stationary state in the joint process $(x_k,\pi_k)$ and note that when $x_k$ is stationary, the Markov chain defined by this process admits an invariant probability measure. The unique ergodicity argument builds on the fact that any two invariant probability measures would have to have the same marginal invariant measure on the state process $x_k$, and this leads to a direct argument in \cite[Lemma B.1]{van2009uniformSPA} to relate unique ergodicity to filter stability. In the context of finite state and measurement spaces, another line of work, adopted in \cite{Kaijser},\cite{kochman2006simple} and \cite{chigansky2010complete}, studies the ergodicity and reachability properties of random matrix products, see \cite{kaijser2011markov} for a countable space setup. In \cite{kochman2006simple} a reachability condition (across all initial priors) through the approximability of a rank-one matrix of unnormalized product of transition matrices, and in \cite{Kaijser} a more restrictive subrectangularity condition, is utilized to establish unique ergodicity. A related argument appears in \cite{Szarek}, in a general context. Finally, \cite{chigansky2010complete} established that the conditions of \cite{kochman2006simple} are tight; with a different sufficiency proof related to filter stability and convex ordering of measures \cite{strassen1965existence}. Going beyond unique ergodicity, geometric ergodicity has recently been studied in \cite{demirci2023geometric}.

\subsection{Statement of main results and contributions}

The paper makes the following contributions.
\begin{itemize}
\item[(i)] We present an approach of defining the state as an infinite-dimensional controlled Markov chain under two classes of policies: a strictly admissible one where the control is a function of the information process or a relaxed one where the control policies satisfy conditional independence between the state and control actions given the information. We provide new sufficient conditions for the existence of optimal control policies which turn out to be stationary in the history variables. In particular, while in the belief-MDP reduction of POMDPs, weak Feller condition requirement imposes total variation continuity on either the system kernel or the measurement kernel, with the approach of this paper only weak continuity of both the transition kernel and the measurement kernel is sufficient together with a stability condition (Theorem \ref{controlledPOMDPE}). 

\item[(ii)] For the discounted cost criterion, more relaxed existence conditions will be presented in Section \ref{discounCostCase}. We first establish the existence of a stationary optimal policy. We then establish near optimality of finite window policies via a direct argument involving near optimality of quantized approximations for MDPs under weak Feller continuity (see Theorem \ref{DiscountExist}), where finite truncations of memory can be viewed as (uniform) quantizations of infinite memory (which is the state definition adopted for the discounted cost criterion by Theorem \ref{discountedMDPC}) with a uniform diameter in each finite window restriction under the product metric. Building on recent results on near-optimality of quantized policies for weak Feller MDPs \cite{SaYuLi15c}, near optimality of finite memory policies follows, complementing and generalizing recent works on the subject \cite{kara2020near,kara2021convergence}. This also facilitates reinforcement learning theoretic methods for POMDPs, which is discussed in further detail in the paper.

\item[(iii)] For the average cost criterion, the paper provides new existence conditions. Furthermore, an approach on how to generate initial priors and beliefs so that optimal performance can be attained is presented: the infinite dimensional formulation presents a natural flexibility in initialization, which can be used to arrive at optimal performance under an absolute continuity condition, see Theorem \ref{FilterInvarianceMeasure2} (which does not seem to be feasible under the standard belief-MDP approach).

\item[(iv)] Via our approach, further results on existence and uniqueness of invariant probability measures for control-free nonlinear filter processes will be established in Section \ref{conFreeCase}. In particular, we will see that unique ergodicity of the measurement process and a measurability condition related to filter stability leads to unique ergodicity (this complements, with an alternative argument, several results in the literature).

%\item[(iv)] An implication of our analysis is that optimal (sliding block) finite memory controllers (with a fixed memory lengh) can be analyzed using the techniques of the paper. We show that under mild conditions, an optimal controller among this class will be stationary. This result is new to our knowledge.
%\item[(iii)] On the traditional method of reducing POMDPs to belief-MDPs, the paper obtains relaxed conditions when compared with the existing results in the literature on the existence of average cost optimal control problems for POMDPs. 
%\item[(vi)] On item (v), we are also able to establish the existence of sample-path optimal controllers under unique ergodicity conditions.
\end{itemize}

Our analysis will begin with the belief-separation approach, and a discussion on its limitations, to be followed with the alternative approach presented in the paper.

\section{Belief-Separation Based Results and Limitations}

As noted earlier, if one wishes to follow the traditional method of reducing a POMDP to a belief-MDP for studying the existence and structure of optimal policies, it would be important to obtain continuity properties of the belief-MDP so that the standard measurable selection theorems, e.g. \cite[Chapter 3]{HernandezLermaMCP} can be invoked to establish the existence of optimal control policies. Building on \cite{KSYWeakFellerSysCont} and \cite{FeKaZg14}, we briefly review the weak Feller property of the kernel defined in (\ref{kernelFilter}) under two different sets of assumptions.

\begin{assumption}\label{TV_channel}
\begin{itemize}
\item[(i)] The transition probability ${\cal T}(\cdot|x,u)$ is weakly continuous in $(x,u)$, i.e., for any $(x_n,u_n)\to (x,u)$, ${\cal T}(\cdot|x_n,u_n)\to {\cal T}(\cdot|x,u)$ weakly.
\item[(ii)] The observation channel $Q(\cdot|x)$ is continuous in total variation, i.e., for any $x_n \to x$, $Q(\cdot|x_n) \rightarrow Q(\cdot|x)$ in total variation.
\end{itemize} 
\end{assumption}

\sy{See \cite[Section 2.1]{KSYContQLearning} for explicit examples on weak continuity and total variation continuity of kernels.}

\begin{assumption}\label{TV_kernel}
The transition probability ${\cal T}(\cdot|x,u)$ is continuous in total variation in $(x,u)$, i.e., for any $(x_n,u_n)\to (x,u)$, ${\cal T}(\cdot|x_n,u_n) \to {\cal T}(\cdot|x,u)$ in total variation.
%\item[(ii)] The observation channel $Q(\cdot|x)$ is independent of the control variable.
\end{assumption}

\begin{theorem}\label{weakFellerBeliefMDP}
\begin{itemize}
\item[(i)] \cite{FeKaZg14}  \label{TV_channel_thm}
Under Assumption \ref{TV_channel}, the transition kernel $\eta(F(\pi,u,Y_{1}) \in \cdot|\pi,u)$ of the filter process is weakly continuous in $(\pi,u)$.
\item[(ii)] \cite{KSYWeakFellerSysCont} Under Assumption \ref{TV_kernel}, the transition kernel $\eta(F(\pi,u,Y_{1}) \in \cdot|\pi,u)$ of the filter process is weakly continuous in $(\pi,u)$.
\end{itemize}
\end{theorem}

We refer the reader to \cite[Theorem 7]{kara2020near}, which builds on \cite{KSYWeakFellerSysCont}, for further refinements with explicit moduli of continuity for the weak Feller property. We refer the reader also to \cite[Theorem 7.1]{FeKaZg14} which establishes weak Feller property under further sets of assumptions. See \cite{feinberg2022markov,feinberg2023equivalent} for additional results on the weak Feller property with a converse theorem statement involving \cite{KSYWeakFellerSysCont} (with regard to both necessity and sufficiency).

\sy{
Since the cost function $c(x,u)$ is continuous and bounded, an application of the generalized dominated convergence theorem\footnote{\begin{theorem}\label{langen}\cite[Theorem 3.5]{Lan81} and \cite[Theorem 3.5]{serfozo1982convergence}Suppose that $\{\mu_n\}_n \subset {\cal P}(\mathbb{X})$ converges weakly to some $\mu$. For a bounded real valued sequence of functions $\{f_n\}_n$ such that $\|f_n\|_\infty <C$ for all $n>0$ with $C<\infty$, if $\lim_{n \to \infty}f_n(x_n)=f(x)$ for all $x_n \to x$, i.e. $f_n$ continuously converges to $f$, then $\lim_{n \to \infty}\int_{\mathbb{X}}f_n(x)\mu_n(dx)=\int_{\mathbb{X}}f(x)\mu(dx)$.
\end{theorem}} implies that $\tilde{c}(z,u)=\int z(dx) c(x,u): {\cal P}(\mathbb{X}) \times \mathbb{U} \to \mathbb{R}$ is also continuous and bounded.}

In the uncontrolled setting, \cite{bhatt2000markov} and \cite{budhiraja2002invariant} have established similar weak continuity conditions (i.e., the weak-Feller property) of the nonlinear filter process in continuous time and discrete time, respectively, where Assumption \ref{TV_channel} is present for an additive noise measurement model.

The convex analytic approach \cite{Manne,Borkar2} is a powerful approach to the optimization of infinite-horizon problems. It is particularly effective in proving results on the optimality of stationary (and possibly randomized stationary) policies, through an infinite-dimensional linear program for constrained optimization problems and infinite horizon average cost optimization problems. 

\begin{assumption}
\label{partial:Implied}
\begin{itemize}
\item[(i)] The $\mathbb{R}_+$-valued one-stage cost function $c$ is bounded and continuous.
\item[(ii)] Assumption \ref{TV_channel} or \ref{TV_kernel} holds (this leads to the belief-MDP $\eta$ to be weakly continuous by Theorem \ref{weakFellerBeliefMDP}).
\item[(iii)] $\mathbb{X}$ and $\mathbb{U}$ are compact.
\end{itemize}
\end{assumption}

%We report the results above as a theorem in the following.

\sy{\begin{theorem}\label{Belief-MDPMethod}
Under Assumption \ref{partial:Implied}, (i) there exists an optimal invariant measure and an associated optimal control policy\footnote{Here, the optimality result may only hold for a restrictive class of initial conditions or initializations.} for the average cost criterion (\ref{expCost}), and (ii) there exists an optimal policy for the discounted cost (\ref{expDiscCost}) criterion.
\end{theorem}
}
\textbf{Proof.} (i) \sy{For the average cost criterion, recall that via (\ref{kernelFilter})-(\ref{weak:eq8}) under the belief-MDP reduction, we are interested in the minimization:
\begin{eqnarray}\label{constOpt1}
\inf_{\gamma \in \Gamma} \limsup_{N \to \infty} {1 \over N} E^{\gamma}_{\mu} [ \sum_{t=0}^{N-1} \tilde{c}(\pi_t,u_t)],
\end{eqnarray}
where $E^{\gamma}_{x_0}[\cdot]$ denotes the expectation over all sample paths with initial distribution $\mu$ so that $x_0 \sim \mu$ under the admissible policy $\gamma$.}
For this criterion, the proof follows from the convex analytic method, studied by Borkar in the weakly continuous case \cite{Borkar2}; for completeness, and its use later in the paper, a proof is presented in Section \ref{proofBelief-MDPMethod}. 

(ii) For the discounted cost criterion, via the belief-MDP (\ref{kernelFilter})-(\ref{weak:eq8}), the result follows from a standard verification and measurable selection analysis on the discounted cost optimality equation \cite[Chapter 3]{HernandezLermaMCP}. \qed

While Theorem \ref{Belief-MDPMethod} presents sufficient findings on the existence of optimal control policies, the approach has a number of limitations: %In the following, we develop an alternative approach with the following benefits to be arrived at:
\begin{itemize}
\item[(i)] [Fully observed MDPs viewed as a special case of POMDPs] Consider the case with the state transition kernel ${\cal T}$ being weakly continuous and where the measurements satisfy $y_t=x_t$, that is, with full state information, in which case the measurement kernel is 
\begin{align}\label{MDPasPOMDP}
Q(dy|x)=\delta_{x}(dy),
\end{align} and thus the system does not satisfy Assumption \ref{TV_channel} or \ref{TV_kernel} (as \sy{the measurement kernel or channel} is only weakly continuous but not total variation continuous). On the other hand, this model, which is a fully observed setup, has been studied extensively and leads to well-known existence (and further, structural and approximation) results. Therefore, we observe that one should be able to relax the conditions further, which we will indeed find to be the case.  
\item[(ii)] For the average cost criterion, Assumption \ref{partial:Implied} provides weaker conditions when compared with the efforts in the literature that typically have utilized the vanishing discount approach:  \cite{Bor00} considers the finite model setup and \cite{borkar2004further} consider restrictive assumptions on the controlled state process, and \cite{feinberg2012average} which in the POMDP case would require yet to be established rates of convergence conditions to invariant measures (needed for uniform boundedness properties of relative discounted cost value functions as the discount parameter approaches unity). We will see that further relaxations can be provided. However, the following limitation, noted in item (iii), is present: 
%Nonetheless, the conditions stated in Theorem \ref{weakFellerBeliefMDP} is restrictive in that the weak-Feller property for the filter process need to be ensured, even though these conditions are weaker than what are presented in the literature, the most relaxed to date being \cite{borkar2004further} which not only requires the total variation continuity of the measurement kernels, but also total variation continuity of the control-kernel. In the following section, we will see that the conditions can be relaxed further.
\item[(iii)] Via the convex analytic method, the obtained policy may not be optimal for all initial states, leading to a constrained notion of optimality. There is then the subtle question of how to select the initial measure, that is the question of achieving the optimal cost under a given initial belief measure and whether the support of an invariant probability measure attracts the belief process from a given initial condition. For an interesting counterexample, see \cite[Example 2.3]{KYSICONPrior}. For belief MDPs, there are few results on unique ergodicity properties, therefore the impact of the initial priors is an open problem under the belief-separated-approach. In our case, the infinite dimensional formulation to be presented leads to a flexibility on how to select the initial prior which we show to lead to global optimality under a mild and testable absolute continuity condition (this approach does not seem to easily carry over to the belief-separated-approach). 
\item[(iv)] Related to (iii), the unique ergodicity problem for the controlled or uncontrolled nonlinear filtering process entails open problems. An alternative approach to what currently exists in the literature is needed, especially for the controlled setup. For the control-free case, building on \cite[Theorem 2]{DiMasiStettner2005ergodicity} and \cite[Prop 2.1]{van2009uniformSPA}, it can be shown that almost sure filter stability in the total variation sense (or weak merging sense, with the arguments tailored to this case) and the uniqueness of an invariant probability measure on the hidden state process leads to unique ergodicity of the filter process. We will see that a complementary condition for the control-free case can be established (unique ergodicity of the measurement process and a measurability condition related to filter stability also leads to unique ergodicity).

%For the control-free case, an ergodicity result can be obtained directly from the ergodicity properties of the process $(X_t,Y_{(-\infty,t]} )$ under an optimal stationary policy, so that the optimal cost can be achieved with an appropriately chosen initial condition/distribution. For the controlled setup, the ergodicity problem is significantly more challenging, but under such a condition the initialization problem can be addressed.% \sy{NOT CLEAR WHAT THE IMPLICATION OF THIS IS..}
\end{itemize}

\sy{Beyond the weak Feller continuity of the kernel $\eta$, if one imposes further regularity, such as Wasserstein regularity, the following holds:
\begin{assumption}\label{main_assumption}
\noindent
\begin{enumerate}
\item \label{compactness}
$(\mathbb{X}, d)$ is a bounded compact metric space 
with diameter $D$ (where $D=\sup_{x,y \in \mathbb{X}} d(x,y)$).
\item \label{totalvar}
The transition probability ${\cal T}(\cdot \mid x, u)$ is 
continuous in total variation in $(x, u)$, i.e., 
for any $\left(x_n, u_n\right) \rightarrow(x, u), 
{\cal T}(\cdot \mid x_n, u_n) \rightarrow 
{\cal T}(\cdot \mid x, u)$ in total variation.
\item \label{regularity}
There exists 
$\alpha \in R^{+}$such that 
$$
\left\|\mathcal{T}(\cdot \mid x, u)-\mathcal{T}\left(\cdot \mid x^{\prime}, u\right)\right\|_{T V} \leq \alpha d(x, x^{\prime})
$$
for every $x,x' \in \mathbb{X}$, $u \in \mathbb{U}$.
\item \label{CostLipschitz}
There exists $K_1 \in \mathbb{R}^+$ such that
\[|c(x,u) - c(x',u)| \leq K_1 d(x,x').\]
for every $x,x' \in \mathbb{X}$, $u \in \mathbb{U}$.
\item The cost function $c$ is bounded and continuous.
%\item \label{dobr_coef}
%$$K_2:=\frac{\alpha D (3-2\delta(Q))}{2}.$$
\end{enumerate}
\end{assumption}

\begin{theorem}\label{ergodicity}\cite{demirci2023average}
    Assume that $\mathbb{X}$ and $\mathbb{Y}$ are Polish spaces. 
    If Assumptions \ref{main_assumption}-\ref{compactness},\ref{regularity} are 
    fulfilled, then we have
    $$
    W_{1}\left(\eta(\cdot \mid z_0, u), \eta\left(\cdot \mid z_0^{\prime},u\right)\right) 
    \leq K_2 W_{1}\left(z_0, z_0^{\prime}\right),$$
with    $K_2:=\frac{\alpha D (3-2\delta(Q))}{2}$
    for all $z_0,z_0' \in \cal{P}(\mathbb{X})$, $u \in \mathbb{U}$.
\end{theorem}

\cite[Theorem 4.37]{SaLiYuSpringer} leads to Lipschitz regularity of the value function under the Wasserstein continuity condition on the kernel:
\begin{theorem}\cite{demirci2023average}
Under Assumption \ref{main_assumption}, for any $\beta \in (0,1)$, there exists an optimal solution to the discounted cost optimality problem with a continuous and bounded value function. Furthermore, with $K_2=\frac{\alpha D (3-2\delta(Q))}{2}$, if $\beta K_2 < 1$ the value function is Lipschitz continuous.
\end{theorem}

The average cost optimality equation (ACOE) plays a crucial role for the analysis and the existence results of MDPs under the infinite horizon average cost optimality criteria. The triplet $(h,\rho^*,\gamma^*)$, where $h,\gamma:\cal{P}(\mathbb{X})\to \mathbb{R}$ are measurable functions and $\rho*\in \mathbb{R}$ is a constant forms the ACOE if 
\begin{align}\label{acoe}
h(z)+\rho^*&=\inf_{u\in\mathbb{U}}\left\{\tilde{c}(z,u) + \int h(z_1)\eta(dz_1|z,u)\right\}\nonumber\\
&=\tilde{c}(z,\gamma^*(z)) + \int h(z_1)\eta(dz_1|z,\gamma^*(z))
\end{align}
for all $z\in\cal{P}(\mathbb{X})$. It is well known that (see e.g. \cite[Theorem 5.2.4]{HernandezLermaMCP}) if (\ref{acoe}) is satisfied with the triplet $(h,\rho^*,\gamma^*)$, and furthermore if $h$ satisfies
\begin{align*}
\sup_{\gamma\in\Gamma}\lim_{t \to \infty}\frac{E_z^\gamma[h(Z_t)] }{t}=0, \quad \forall z\in\cal{P}(\mathbb{X})
\end{align*}
then $\gamma^*$ is an optimal policy for the POMDP under the infinite horizon average cost optimality criteria, and 
\begin{align*}
J_{\infty}^*(z)=\inf_{\gamma\in\Gamma}J(z,\gamma)=\rho^* \quad \forall z\in \cal{P}(\mathbb{X}).
\end{align*}
%We will refer to the function $h$, the relative value function, for the rest of the paper. Note that there may not be a unique relative value function $h$ that satisfies %the ACOE, however, any $h$ that satisfies the ACOE can be used for optimality analysis.
%We thus state the following main theorem.
\begin{theorem}\label{mainEmre} \cite{demirci2023average}
Under Assumption \ref{main_assumption}, with $K_2=\frac{\alpha D (3-2\delta(Q))}{2} < 1$, 
    a solution to the average cost optimality 
    equation (ACOE) exists. 
    This leads to the existence of an optimal 
    control policy, and optimal cost is constant for 
    every initial state.
\end{theorem}
In the following, these will be relaxed significantly. 
}

\section{Another Look at POMDPs: Infinite Dimensional Markov Decision Process Formulation}

To facilitate the approach presented in this paper, we first assume that the measurements and the control actions have been taking place since $- \infty$. Later on we will motivate and justify this approach. 

%But is this realizable? That is, at time $0$, there is no information at all! Thus the initial belief is just the marginal probability at time zero on $X_0$ in the Bayesian update equation. We can expand the measurement space to include a special symbol $N$ which takes place with zero measure. In the end, this would just be an initial condition and the mathematics of the analysis would evolve regardless of where we start from. In particular the convex analytic method does not depend on the initial condition as long as every sequence has a converging subsequence, which is the case in our setup. So, this won't be a problem.}

Let $\mathbb{Y}^{\mathbb{Z}_-}$ be the one-sided product space consisting of elements of the form 
\[y_{(-\infty,0]} = \{\cdots, y_n, \cdots, y_{-2}, y_{-1}, y_0\}\]
 with $y_k \in \mathbb{Y}$. We endow $\mathbb{Y}^{\mathbb{Z}_-}$ with the product topology; this makes $\mathbb{Y}^{\mathbb{Z}_-}$ a metric space, which is complete and separable. Likewise, we will view $U_{(-\infty,-1]}$ also as a $\mathbb{U}^{\mathbb{Z}_-}$-valued random variable. \sy{On $\mathbb{Y}^{\mathbb{Z}_-} \times \mathbb{U}^{\mathbb{Z}_-}$, we define the following product metric: For $\bar{y}_{(-\infty,0]},\bar{u}_{(-\infty,-1]} \in \mathbb{Y}^{\mathbb{Z}_-} \times \mathbb{U}^{\mathbb{Z}_-}$,
\begin{align}
&\bar{d}(\bar{y}_{(-\infty,0]},\bar{u}_{(-\infty,-1]},\tilde{y}_{(-\infty,0]},\tilde{u}_{(-\infty,-1]})  \nonumber \\
& \qquad := \sum_{k=0}^{\infty} 2^{-k} \frac{d_{\mathbb{Y}}(\bar{y}_{-k},\tilde{y}_{-k})}{1+d_{\mathbb{Y}}(\bar{y}_{-k},\tilde{y}_{-k})} + \sum_{m=1}^{\infty} 2^{-(m-1)} \frac{d_{\mathbb{U}}(\bar{u}_{-m},\tilde{u}_{-m})}{1+d_{\mathbb{U}}(\bar{u}_{-m},\tilde{u}_{-m})},
\end{align}
 }
A related view, on using the infinite past as a Markov process, was presented in \cite{yukselSICON2017} to establish stochastic stability of control-free dynamical systems driven by noise processes which are not necessarily independent and identically distributed, but which are stationary. A further related approach was adopted in \cite{albertini2001small}; see also \cite{HairerErgodic} (as well as \cite{bismut1982partially} and \cite{YukselWitsenStandardArXiv} as discussed earlier).

\sy{
We will construct two controlled Markov models: The first one will primarily be useful for the discounted cost criterion and in some ways more natural; the second one will lead to a more relaxed formulation with stronger results for the average cost criterion (while the former is also applicable for the average cost criterion under more restrictive conditions).  
}

 \subsection{Infinite-product information memory as controlled state: Strong sense admissible policies}
  In our first formulation, which will be utilized in the discounted infinite horizon optimality analysis the {\it state variable}
 will be defined as
  \[S_k:= (Y_{(-\infty,k]}, U_{(-\infty,k-1]}),\] 
which is the information available to the controller at time $k$, provided that the infinite past history is available. We have the following immediate result:

 \begin{lemma}\label{MDPwithSk}
The process $(S_k, U_k)$ is a controlled Markov chain where $S_k$ is $\mathbb{Y}^{\mathbb{Z}_-} \times \mathbb{U}^{\mathbb{Z}_-}$-valued.
\end{lemma}
 
\textbf{Proof.} 
Let $A^y = \prod_{k=-\infty}^{0} A^y_k$ be an open cylinder set in $\mathbb{Y}^{\mathbb{Z}_-}$ and $A^u = \prod_{k=-\infty}^{0} A_k^u$ be an open cylinder set in $\mathbb{U}^{\mathbb{Z}_-}$: We have that $P$-almost surely
\begin{eqnarray}\label{MarkovKernel1}
&& P\bigg(S_{t+1} \in (A^y \times A^u) | S_t=s, U_t=u, S_{[0,t-1]}=s_{[0,t-1]}, U_{[0,t-1]}=u_{[0,t-1]}\bigg) \nonumber \\
&& = \int P(Y_{t+1} \in A^y_0 | y_{(-\infty,t]}, u_{(-\infty,t]}) P( Y_{(-\infty,t]} \in \prod_{k=-\infty}^{-1} A^y_k  | y_{(-\infty,t]}) P(U_{(-\infty,t]} \in \prod_{k=-\infty}^{-1} A^u_k  | u_{(-\infty,t]})   \nonumber \\
%&& = \int P(y_{t+1} \in A^y_0 | y_{(-\infty,t]}, u_{(-\infty,t]}) P( y_{(-\infty,t]} \in \prod_{k=-\infty}^{-1} A^y_k  | y_{(-\infty,t]}) P(u_{(-\infty,t]} \in \prod_{k=-\infty}^{-1} A^u_k  | u_{(-\infty,t]})   \nonumber \\
&& P\bigg(S_{t+1} \in (A^y \times A^u) | S_t=s, U_t=u \bigg) \nonumber \\
&&=: \mathbb{P}^d(A^y \times A^u | S_t=s, U_t=u)
\end{eqnarray}
\qed

We will use this construction in the following section; together with near optimality results of quantized MDPs with respect to the product topology, we will be able to obtain rigorous and very general near optimality results of sliding finite window based control policies.

 \subsection{A relaxation useful for average cost control}\label{sectionMDPAC}
We will also study an alternative construction which will prove to be useful for the average cost criterion, towards a relaxed set of admissible policies.

%\sy{Compared with \cite[Section 3.2]{YukselWitsenStandardArXiv}, we don't require conditional independence from the entire past; just the state. This is because the analysis here is tailored towards POMDPs, and not decentralized stochastic control where a state variable is more general. If we also make it independent from the entire past, this will make the analysis more tedious for average cost problems, though this is possible. But because we can establish that an optimal policy will be deterministic using our paper with Ari, it will be all right; as deterministic policies will be extremal and attain optimality.}

%\sy{if we replace $\omega_0$ with $x_0, w_{[0,\infty]}$, it will be a special case of our SIAM setup. It is more convenient here to use the more relaxed version with the state. Randomization is well-defined anyway, so no problem; just independent from the entire past!!! Use this, it will be a good contribution...Read Bismut again also.}

\begin{lemma}\label{CMCZ_k}
With $Z_k= (Y_{(-\infty,k]}, U_{(-\infty,k-1]},X_k)$, the pair $(Z_k, U_k)$ is a controlled Markov chain where $Z_k$ is $\mathbb{Y}^{\mathbb{Z}_-} \times \mathbb{U}^{\mathbb{Z}_-} \times \mathbb{X}$ valued.
\end{lemma}

\textbf{Proof.} Let $A = \prod_{k=-\infty}^{0} A^y_k$ be an open cylinder set in $\mathbb{Y}^{\mathbb{Z}_-}$ and $A^u = \prod_{k=-\infty}^{0} A_k^u$ be an open cylinder set in $\mathbb{U}^{\mathbb{Z}_-}$ and let $B \in \mathbb{X}$ be a Borel set. We have that $P$-almost surely
\begin{eqnarray}\label{MarkovKernel}
&& P\bigg(Z_{t+1} \in (A^y \times A^u \times B) | Z_t=z, U_t=u, Z_{[0,t-1]}=z_{[0,t-1]}, U_{[0,t-1]}=u_{[0,t-1]}\bigg) \nonumber \\
&& = \int Q(Y_{t+1} \in A^y_0 | x_{t+1}) {\cal T}(X_{t+1} \in B | x_t=x, u_t=u) \nonumber \\
&& \quad \quad \quad \quad \quad \quad \times P( Y_{(-\infty,t]} \in \prod_{k=-\infty}^{-1} A^y_k  | y_{(-\infty,t]}) P(U_{(-\infty,t]} \in \prod_{k=-\infty}^{-1} A^u_k  | u_{(-\infty,t]})   \nonumber \\
&& = \int Q(Y_{t+1} \in A^y_0 | x_{t+1}) {\cal T}(X_{t+1} \in B | x_t=x, u_t=u) 1_{\{ y_{(-\infty,t]} \in \prod_{k=-\infty}^{0} A^y_k, u_{(-\infty,t]} \in \prod_{k=-\infty}^{0} A^u_k\} } \nonumber \\ \label{MarkovKerne2} \\
&& =: \mathbb{P}^a(A^y \times A^u \times B | Z_t=z, U_t=u)
\end{eqnarray}
where for a Borel $C$, $Q(y \in C|x) := P(G(X_0,V_0) \in C | X_0=x)$ by (\ref{updateEq2}) and ${\cal T}$ is the transition kernel defined by (\ref{updateEq1}).
\qed

In the above, we define $\mathbb{P}^a$ to be the transition kernel for this Markov chain. 

\begin{definition}[Weak-Feller Condition for $\mathbb{P}^a$]
$\mathbb{P}^a(dz_{t+1}  |z_t, u_t)$ is weak-Feller if for every continuous and bounded $f \in C_b(\mathbb{Y}^{\mathbb{Z}_-} \times \mathbb{U}^{\mathbb{Z}_-} \times \mathbb{X})$:
\[\int f(z_{t+1}) \mathbb{P}^a(dz_{t+1} | z_t=z, u_t=u)\]
is continuous in $(z,u)$.
\end{definition}

The following will be invoked often.
\begin{assumption}\label{contCond}
$Q(\cdot | x)$ is weakly continuous in $x$, and ${\cal T}(\cdot | x,u)$ is weakly continuous in $x,u$.
\end{assumption}

\begin{lemma}\label{weakFellerP}
Under Assumption \ref{contCond}, $\mathbb{P}^a$ is weak-Feller.
\end{lemma}

\textbf{Proof.}
First note that by (\ref{MarkovKernel}) we have
\[\mathbb{P}^a(dz_{t+1} | z_t, u_t) = Q(dy_{t+1}|x_{t+1}) {\cal T}(dx_{t+1} | x_t, u_t) P(dy_{(-\infty,t]} | y_{(-\infty,t]} ) P(du_{(-\infty,t]} | u_{(-\infty,t]} ) \]
where $P(y_{(-\infty,t]} \in \cdot | y_{(-\infty,t]}) = \delta_{y_{(-\infty,t]}}(\cdot)$ is the Dirac measure (as in (\ref{MarkovKerne2})), with the same holding for $P(du_{(-\infty,t]} | u_{(-\infty,t]})$.

We wish to show that the following is continuous in $x_{t}, y_{(-\infty,t]},u_{(-\infty,t-1]}, u_t$, for continuous and bounded $f$:
\begin{eqnarray}
&&\int f\bigg(x_{t+1}, y_{(-\infty,t+1]},u_{(-\infty,t]} \bigg) \mathbb{P}^a(dx_{t+1}, dy_{(-\infty,t]+1},du_{(-\infty,t-1} | x_{t}, y_{(-\infty,t]},u_{(-\infty,t-1}, u_t) \nonumber \\
&&= \int f\bigg(x_{t+1}, y_{(-\infty,t+1]},u_{(-\infty,t]} \bigg)  Q(dy_{t+1}|x_{t+1}) {\cal T}(dx_{t+1} | x_t, u_t) \nonumber \\
&& \qquad \qquad \qquad \times P(dy_{(-\infty,t]} | y_{(-\infty,t]} ) P(du_{(-\infty,t]} | u_{(-\infty,t]} ) \nonumber \\
&&= \int \int_{x_{t+1}} \bigg( \bigg( \int_{y_{t+1}}  f(x_{t+1}, y_{(-\infty,t+1]},u_{(-\infty,t]})  Q(dy_{t+1}|x_{t+1}) \bigg) {\cal T}(dx_{t+1} | x_t, u_t) \bigg) \nonumber \\
&& \qquad \qquad \qquad \times P(dy_{(-\infty,t]} | y_{(-\infty,t]} ) P(du_{(-\infty,t]} | u_{(-\infty,t]} ) \nonumber 
\end{eqnarray}
The expression
\[ g(x_{t+1},y_{(-\infty,t+1]},u_{(-\infty,t}) := \bigg( \int_{y_{t+1}}  f(x_{t+1}, y_{(-\infty,t+1]},u_{(-\infty,t})  Q(dy_{t+1}|x_{t+1}) \bigg) \]
is continuous by the generalized dominated convergence theorem (Theorem \ref{langen}) for varying measures under {\it continuous-convergence} \cite[Thm. 3.5]{Lan81} (see \cite[Thm. 3.3]{serfozo1982convergence} for a more restricted version). Since $g$ is continuous in $x_{t+1}$, 
\[h(x_t, u_t, y_{(-\infty,t]},u_{(-\infty,t]}) := \int_{x_{t+1}} \bigg( g(x_{t+1},y_{(-\infty,t]},u_{(-\infty,t]}) {\cal T}(dx_{t+1} | x_t, u_t) \bigg)\]
is also continuous in its parameters. Finally, again by \cite[Thm. 3.5]{Lan81},
\[\int h(x_t, u_t, y_{(-\infty,t]},u_{(-\infty,t]}) P(dy_{(-\infty,t]} | y_{(-\infty,t]} ) P(du_{(-\infty,t]} | u_{(-\infty,t]} )\]
is continuous in $x_{t}, y_{(-\infty,t]},u_{(-\infty,t-1]}, u_t$.
%\[\int \int_{x_{t+1}} \bigg( g(x_{t+1},y_{(-\infty,t]},u_{(-\infty,t]}) P(dx_{t+1} | x_t, u_t) \bigg) P(dy_{(-\infty,t]} | y_{(-\infty,t]} ) P(du_{(-\infty,t]} | u_{(-\infty,t]} ) \nonumber \\
%\]
\qed

The above completes the description of the controlled Markov model, that we will utilize for the average cost problem. 

%Note that we have not yet utilized the relaxed control policy framework (i.e., how the actions are generated). As noted earlier, for the average cost criterion, we will later view the policy space as the set of all joint measures that satisfies 
%\begin{eqnarray}\label{ConditionallyHistoryIndependentPolicies}
%X_n \leftrightarrow (U_{(-\infty,n-1]}),Y_{(-\infty,n]}) \leftrightarrow U_n
%\end{eqnarray}
%For every fixed random variable $(U_{(-\infty,n-1]}),Y_{(-\infty,n]}, X_n)$, by viewing this conditional independence as a special instance of \cite[Theorem 5.2(i)]{saldiyukselGeoInfoStructure}, it follows that this policy space is weakly compact at every time instant with a fixed measure on $(U_{(-\infty,n-1]}),Y_{(-\infty,n]}, X_n)$. 

\subsection{Key technical assumptions on continuity and stability}

For some of the results to be presented, we will invoke some (but not all) of the following assumptions on stationarity, stability, and continuity. %We will see, in Section \ref{sufficientCondKey}, that these are related to filter stability \cite{chigansky2009intrinsic}. 

\begin{assumption}\label{contConditinal}[A Continuity Condition]
The stochastic kernel $P(dx_k | y_{(-\infty,k]},u_{(-\infty,k-1]})$ is weakly continuous, that is, for every continuous and bounded $g: \mathbb{X} \to \mathbb{R}$,
\[\int g(x) P(X_k \in dx | y_{(-\infty,k]},u_{(-\infty,k-1]})\] is continuous in $(y_{(-\infty,k]},u_{(-\infty,k-1]})$ (under the product topology on $\mathbb{Y}^{\mathbb{Z}_-} \times \mathbb{U}^{\mathbb{Z}_-}$).
\end{assumption}

\begin{assumption}\label{assumptionMC}[A Filter Stability Condition]
$\pi_t$ is $\sigma(I_t)$ measurable where $I_t = (Y_{( -\infty,t]},U_{( -\infty,t-1]})$. In particular, independent of the control policy, there exists a measurable $g$ such that $g(y_{(-\infty,t]},u_{(-\infty,t-1]})(A) = \pi_t(A)$ for every Borel $A$. 
\end{assumption}

\begin{assumption}\label{ForgetCondition}[A Stationarity / Stability Condition]
For every Borel $B$, and with $y = (\cdots, y_{-2}, y_{-1}, y_0)$, $u = (\cdots, u_{-2}, u_{-1}, u_0)$, under any control policy, for almost every $(y, u)$, we have that
\[P(X_t \in B| Y_{(-\infty,t]}=y, U_{(-\infty,t-1]}=u) = P(X_{t+1} \in B| Y_{(-\infty,t+1]}=y, U_{(-\infty,t]} = u)\]
\end{assumption}

\sy{Further discussion on these assumptions is given in Appendix \ref{OnAssumptions}, where we point out that some of the conditions above are related to {\it filter stability} \cite{chigansky2009intrinsic}. Notably, conditions which lead to contraction under the Hilbert projective metric \cite{le2004stability} (which are applicable also in the controlled setup under mixing conditions \cite{demirciRefined2023}, uniform over all control actions), lead to each of the conditions above. }

\section{Discounted Cost Criterion: Existence Results and Near Optimality of Finite Window/Memory Policies}\label{discounCostCase}

In this section, we study the discounted cost setup. We will see a particular utility in this criterion: near optimality of finite window policies. 

%\subsection{Controlled-Markov construction with infinite-memory state space for the discounted cost criterion}
Recall Lemma \ref{MDPwithSk} where the controlled state variable is given with $S_k= (Y_{(-\infty,k]}, U_{(-\infty,k-1]})$. 

Given the construction and stability properties, for the discounted cost criterion, towards arriving at approximate optimality results using finite window policies, we will present an alternative controlled Markov model. The proof follows from identical arguments presented in Section \ref{sectionMDPAC}, but here we apriori impose Assumption \ref{contConditinal} and Assumption \ref{assumptionMC}. 

\begin{theorem}\label{discountedMDPC}
\begin{itemize}
%\item[(i)] With \[S_k= (Y_{(-\infty,k]}, U_{(-\infty,k-1]}),\] the pair $(S_k, U_k)$ is a controlled Markov chain where $S_k$ is $\mathbb{Y}^{\mathbb{Z}_-} \times \mathbb{U}^{\mathbb{Z}_-}$ valued.
\item[(i)] Under Assumption \ref{assumptionMC}, we have that for all $u \in \mathbb{U}$
\begin{align}\label{barcInfProd}
E\bigg[c(X_k,u) \bigg| (Y_{(-\infty,k]}, U_{(-\infty,k-1]}=s) \bigg] =: \bar{c}(s,u),
\end{align}
for some measurable $\bar{c}: \mathbb{Y}^{\mathbb{Z}_-} \times \mathbb{U}^{\mathbb{Z}_-} \times \mathbb{U} \to \mathbb{R}$.
\item[(ii)] Under Assumption \ref{contConditinal}, $\bar{c}$ is continuous and bounded, with $c$ bounded continuous.
\item[(iii)] \sy{The kernel $\mathbb{P}^d(S_{t+1} \in \cdot | S_t=s, U_t=u)$ is weak Feller under Assumptions \ref{contCond} and \ref{contConditinal}.}
\end{itemize}
\end{theorem}
\sy{
\textbf{Proof.} (i) follows from Assumption \ref{assumptionMC}. To see this, note that
\[E\bigg[c(X_k,u) \bigg| (Y_{(-\infty,k]}, U_{(-\infty,k-1]}=s) \bigg] = \tilde{c}(\pi_k,u),\]
Now, by Theorem \ref{langen},
\[\tilde{c}(\pi,u) = \int \pi(dx) c(x,u)\] 
is a continuous function from ${\cal P}(\mathbb{X}) \times \mathbb{U} \to \mathbb{R}$. By representing this continuous function as the  pointwise limit of simple functions $\sum_{k=1}^N \alpha^N_k 1_{\{\pi \in A^N_k, u \in B^N_k\}}$ for Borel sets $\{A^N_k, B^N_k; k=1,\cdots,N\},  N \in \mathbb{N} \}$ and $\{\{\alpha^N_k, k=1, \cdots, N\}; N \in \mathbb{N} \}$, by Assumption \ref{assumptionMC} with $\pi = g(s)$ such that $\sum_{k=1}^N \alpha^N_k 1_{s \in g^{-1}(A^N_k), u \in B^N_k}$ is Borel in the space $\mathbb{Y}^{\mathbb{Z}_-} \times \mathbb{U}^{\mathbb{Z}_-} \times \mathbb{U}$, it follows that the pointwise limit of these measurable functions defining $\bar{c}: (s,u) \mapsto \tilde{c}(g(s),u)$ is also measurable as a map from $\mathbb{Y}^{\mathbb{Z}_-} \times \mathbb{U}^{\mathbb{Z}_-} \times \mathbb{U} \to \mathbb{R}$. 
%Via \cite[Theorem 2.1]{DubinsFreedman},
   (ii) follows from Assumption \ref{contConditinal} by definition. (iii) follows by considering $P(Y_t \in dy | S_t=s, U_t=u) = \int Q(dy|x) P(X_t \in dx| S_t=s, U_t=u)$, and the weak continuity of each of the kernels in the integration and generalized weak convergence \cite[Thm. 3.5]{Lan81} as used earlier in the proof of Lemma \ref{weakFellerP}; 
\qed
}

%\sy{NOTE: The above then also implies the convex analytic method to be applicable since we have all the conditions under weak Feller continuity. SHOULD WE NOTE THIS? PERHAPS A RESTRUCTURING IS NEEDED.}

As a result, under Assumptions \ref{contConditinal} and \ref{assumptionMC}, we have an equivalent controlled Markov model with cost $\bar{c}: (\mathbb{Y}^{\mathbb{Z}_-} \times \mathbb{U}^{\mathbb{Z}_-}) \times \mathbb{U} \to \mathbb{R}_+$, $\mathbb{Y}^{\mathbb{Z}_-} \times \mathbb{U}^{\mathbb{Z}_-}$-valued state $s_k$, and transition kernel  $\mathbb{P}^d(S_{t+1} \in \cdot | s_t=s, U_t=u)$.

% (which has recently been studied in \cite{kara2020near} and \cite{kara2021convergence}). 

For the discounted criterion, we will restrict our control policies to those that are strict-sense admissible. We could also consider the relaxed framework, however the analysis for the discounted criterion setup will be seen to be more direct. Compared with the average cost criterion, we are able to utilize a verification theorem directly following the discounted cost optimality equation (DCOE). 
%Note that arriving at the counterpart of this equation in average-cost, known as the average cost optimality equation or inequality (ACOE/I) \cite{survey}, often requires contraction or value iteration based methods (see e.g.
%\cite{Veg03,hernandez2012adaptive,ABor-19}), or the vanishing discount method (see e.g. \cite{survey,feinberg2012average,HernandezLermaMCP,hernandezlasserre1999further,costa2012average,GoHe95,yu2020average} which have presented various conditions and relaxations). Efforts under this method typically require some recurrence/ergodicity/Dobrushin type geometric or subgeometric convergence conditions,
%which may be too strong for a large class of applications (e.g., for belief-MDP reduction of Partially Observable Markov Decision Processes). This explains the different control policy class we have adopted for the average cost setup, as well as the adoption of the convex analytic method. 

Let, as in Theorem \ref{discountedMDPC}, $s_k=(y_{(-\infty,k]},\tilde{u}_{(-\infty,k-1]}) \in \mathbb{Y}^{\mathbb{Z}_-} \times \mathbb{U}^{\mathbb{Z}_-}$. Now, under strict-sense policies, with the assumptions on weak-Feller continuity and bounded continuous cost function, if there exists a solution to
\[J_{\beta}(s) = \min_{u_k \in \mathbb{U}} E\bigg[c(x_k,u_k) + \beta J_{\beta}(S_{k+1}) | S_k=s, U_k=u\bigg], \]
then we can declare this policy to be optimal using a standard verification argument \cite{HernandezLermaMCP}. As before, we assume that $\mathbb{X}, \mathbb{Y}, \mathbb{U}$ are compact.

\begin{theorem}\label{DiscountExist}
An optimal solution exists under Assumptions \ref{contCond} and \ref{contConditinal} and Assumption \ref{assumptionMC}.
\end{theorem}

\textbf{Proof.} Recall first that 
\[ E[(c(X_k,U_k)) | S_k=s, U_k=u] \]
is a measurable function of $(s,u)$ under Assumption \ref{assumptionMC}, call this $\bar{c}(s,u)$. We then have
\begin{align}\label{DCOEInfP}
J_{\beta}(s) = \min_{u \in \mathbb{U}}  \bigg(E[\bar{c}(S_k,U_k)+ \beta J_{\beta}(S_{k+1}) | S_k=s, U_k=u] \bigg) 
\end{align}
Then, if we have weak-Feller continuity of the controlled kernel (by Assumption \ref{contCond}) and if we have the continuity of $f$ (by Assumption \ref{contConditinal}), by measurable selection theorems \cite{HernandezLermaMCP}, there exists a solution to the discounted cost optimality equation.\qed

\sy{While we have an existence result, there is the operational question that in a given application, (i) we do not have the past data and we only have that $X_0 \sim \mu$ for some $\mu$ and the measurements following $t=0$; and that (ii) computationally, the problem is still demanding with no apriori benefit of this approach. 

In the following, we address both of these issues:

On item (i); now that we have an optimal solution given $s_t$, with (\ref{kernelFilter}) and (\ref{weak:eq8}) and the belief-MDP in $({\cal P}(\mathbb{X}),\mathbb{U},\eta,\tilde{c})$, by a measurable selection theorem of Blackwell and Ryll-Nardzewski \cite{Blackwell3} (see also p. 255 of \cite{dynkin1979controlled}), we first note that for any finite horizon problem, for every Markov policy which uses $s_t$, there exists another Markov policy which uses $\pi_t$ which is at least as good as the policy using $s_t$ \cite{Blackwell2}. Thus, this implies that for any finite horizon problem we can use $\pi_t$ as a sufficient statistic, and an equivalence class among the history process $s_t$ based on their induced $\pi_t$ values can be used to execute an optimal policy. This argument first applies for finite horizon problems, and by value iteration defining a contraction on continuous and bounded functions on both $s_t$ and $\pi_t$, by a standard verification theorem, the result also applies for infinite horizon problems with policies now being stationary. An optimal policy will use $\pi_t$ as a sufficient statistic, and therefore any realization of $s_0$ compatible with $X_0 \sim \mu$ will be optimal.

More operationally, on (ii) as a consequential application, we note the following. We say that, for some $N \in \mathbb{N}$, a control policy is an $N$-memory policy if for $t > N$,
\[U_t = \gamma_t(Y_{[t-N+1,t]},U_{[t-N+1,t-1]}),\] 
and for $0 \leq t \leq N$, 
\[U_t = \gamma_t(Y_{[0,t]},U_{[0,t-1]}).\]

\begin{theorem}\label{finiteWApp}
Under the conditions of Theorem \ref{DiscountExist}, the set of $N$-memory policies are asymptotically optimal; i.e., for every $\epsilon > 0$, there exists $N$ so that $N$ window policies are $\epsilon$-optimal. 
\end{theorem}
}
Before the prove the theorem, let us note that the near-optimality is operationally justified in two scenarios. If the costs start becoming active from time $N$ onwards, or the policies in the first $N$ stages are designed via a policy-iteration argument; see \cite[p. 16]{kara2020near}.

\textbf{Proof.} Finite window truncation of $z$ can be viewed as a uniform binning/quantization (even though the number of bins may not be countable) operation applied to the infinite-dimensional state vector $z$ under the product topology. Indeed, define (the binning/quantization map)
\[\rho: \mathbb{Y}^{\mathbb{Z}_-} \times \mathbb{U}^{\mathbb{Z}_-} \to \mathbb{Y}^{N+1} \times \mathbb{U}^{N},\]
with
\[\rho(y_{(-\infty,0]},\tilde{u}_{(-\infty,-1]}) = (y_{[-N,0]},u_{[-N,-1]})\]
In this case, for any $(\bar{y}_{(-\infty,0]},\bar{u}_{(-\infty,-1]})$ and $(\tilde{y}_{(-\infty,0]},\tilde{u}_{(-\infty,-1]})$ with 
\[\rho(\bar{y}_{(-\infty,0]},\bar{u}_{(-\infty,-1]})=\rho(\tilde{y}_{(-\infty,0]},\tilde{u}_{(-\infty,-1]}),\] and $\bar{d}$ being the product metric defined by
\begin{align}
&\bar{d}(\bar{y}_{(-\infty,0]},\bar{u}_{(-\infty,-1]},\tilde{y}_{(-\infty,0]},\tilde{u}_{(-\infty,-1]})  \nonumber \\
& \qquad := \sum_{k=0}^{\infty} 2^{-k} \frac{d_{\mathbb{Y}}(\bar{y}_{-k},\tilde{y}_{-k})}{1+d_{\mathbb{Y}}(\bar{y}_{-k},\tilde{y}_{-k})} + \sum_{m=1}^{\infty} 2^{-(m-1)} \frac{d_{\mathbb{U}}(\bar{u}_{-m},\tilde{u}_{-m})}{1+d_{\mathbb{U}}(\bar{u}_{-m},\tilde{u}_{-m})},
\end{align}
where $d_{\mathbb{Y}}$ and $d_{\mathbb{U}}$ are the metrics on $\mathbb{Y}$ and $\mathbb{U}$, respectively;
we have that 
\begin{align}\label{windowSizeBin}
\bar{d}(\bar{y}_{(-\infty,0]},\bar{u}_{(-\infty,-1]},\tilde{y}_{(-\infty,0]},\tilde{u}_{(-\infty,-1]})  \leq 3 \frac{1}{2^N}.
\end{align}
\sy{Since the constructed kernel $\mathbb{P}^d$ is weakly continuous}, the cost function is bounded continuous, and the state space is compact, the proof is then a corollary of \cite[Theorem 4.27]{SaLiYuSpringer} (see also \cite[Theorem 11]{KaraYuksel2021Chapter}) by viewing the finite window truncation as a quantization (with a uniformly bounded radius for each bin) of the state space under the product topology. 

In particular, one can construct an approximate MDP model with state space $\mathbb{Y}^{N} \times \mathbb{U}^{N-1}$ and action space $\mathbb{U}$ and transition kernel $\mathbb{P}^N$ whose solution can be extended, via the quantization rule $\rho$, to the whole $\mathbb{Y}^{\mathbb{Z}_-} \times \mathbb{U}^{\mathbb{Z}_-}$ and which will be near optimal for the original problem \cite[Theorem 4.27]{SaLiYuSpringer}.
\qed

%\sy{check $s_k$ vs. $z_k$ (which contains $x_k$ in it, in the earlier statements; to be corrected.)}

{\bf Algorithmic and Numerical Implications.} Note that a finite memory truncation leads to a uniform quantization error, therefore the mathematical analysis adopted in this paper is particularly suitable for such a problem in view of recent near optimality results for quantized approximations to weakly continuous kernels. We note that Theorem \ref{finiteWApp} has a rather direct relation with \cite{kara2020near} and \cite{kara2021convergence}, where near optimality of finite window policies were established via alternative and slightly more tedious methods, and more restrictive conditions (though with rates of convergence properties, which we do not present). In our current paper, the filter stability condition manifests itself via Assumption \ref{assumptionMC} {\it with no additional assumptions other than the weak Feller property}. 

The result above thus presents further sufficient conditions on the applicability of re-inforcement learning methods presented in \cite{kara2021convergence} for finite memory near-optimal control. In particular, if one applies an independently randomized exploration control policy for $u_t$ (for learning) in the $Q$-learning algorithm presented in \cite[Section 4]{kara2021convergence}, one arrives at a control-free system under which the unique ergodicity of the measurement process is satisfied under mild conditions so that the conditions in \cite[Theorem 4.1]{kara2021convergence} (and in this particular fully observed interpretation, \cite[Corollary 3.3]{KSYContQLearning}) is applicable for arriving at near-optimal control policies which use only a finite window of recent measurement and past actions.

\sy{Since finite memory truncation is viewed as quantization, under further regularity conditions on the transition kernels, \cite[Theorem 5.2]{SaYuLi15c} or \cite[Theorem 2.5]{KSYContQLearning} can be used to arrive at rates of convergence with respect to the window size. We also note here that another positive attribute for the discounted criterion is that under weak Feller continuity, both the value functions and optimal policies are robust to model approximations \cite[Theorem 4.4]{kara2020robustness}. }%Therefore, small perturbations in state realizations do not impact the expected costs %significantly. 

\sy{
Going beyond the weak Feller property, we have the following.

\begin{assumption}\label{Wasser1}
We assume that for some $\bar{K}_1, \bar{K}_2$:
\begin{itemize}
\item[(a)] With $\bar{c}$ defined in (\ref{barcInfProd}),
 $|\bar{c}(s,u) - \bar{c}(s',u)| \leq \bar{K}_1 \bar{d}(s,s')$; that is, $\bar{c}(\cdot,u)$ is $\bar{K}_1$-Lipschitz.
\item[(b)] $W_1(\mathbb{P}^d(ds_1 | s_0=s,u_0=u),\mathbb{P}^d(ds_1 | s_0=s',u_0=u)) \leq \bar{K}_2\bar{d}(s,s')$.
\end{itemize}
\end{assumption}

\begin{theorem}\label{LipschitzValueFInf}
Suppose that Assumptions \ref{contCond}, \ref{contConditinal} and \ref{Wasser1} hold. Then, the solution to (\ref{DCOEInfP}), which is the optimal value function, is Lipschitz with coefficient $\bar{K} = \frac{K_1}{1 - \beta \bar{K}_2}$.
\end{theorem}

\textbf{Proof.} Under Assumptions \ref{contCond} and \ref{contConditinal}. The kernel $\mathbb{P}^d(S_{t+1} \in \cdot | S_t=s, U_t=u)$ is weak Feller; see Theorem  \ref{discountedMDPC}(iii). The result then follows from \cite[Theorem 4.37]{SaLiYuSpringer}.
\qed

\cite[Theorem 6]{KSYContQLearning} then implies the following:
\[ \sup _{s \in  \mathbb{Y}^{\mathbb{Z}_-} \times \mathbb{U}^{\mathbb{Z}_-}}\left|J_\beta\left(s, \hat{\gamma}\right)-J_\beta^*\left(s\right)\right| \leq \frac{2 \bar{K}_1}{(1-\beta)^2(1-\beta \bar{K}_2)} \bar{L}(N), \]
where $\hat{\gamma}$ denotes the optimal policy of 
the finite-state approximate model extended to 
the state space and $\bar{L}(N)$ denotes the diameter of the quantization bins, in this case upper bounded by (\ref{windowSizeBin}). Therefore, the quantized model, which gives a finite window approximate policy, is near optimal with an explicit rate of convergence as $N \to \infty$. 

%When $\mathbb{X}$ is not compact, note that we may still need to verify (\ref{OptDisRestricted}) to claim optimality. Here the cost is bounded, so OK.
We leave explicit conditions for Assumption \ref{Wasser1} for future work; however, an immediate sufficient condition is via the projective Hilbert metric \cite{le2004stability} which can be shown to be applicable also in the controlled setup under uniform mixing conditions \cite{demirciRefined2023}, see Appendix \ref{OnAssumptions}.
%We note that the benefit of this additional structure is that one can establish uniform error bounds on the performance of finite memory policies: \cite[Theorem 6]{KSYContQLearning} implies that the loss is linear in the approximation error. 

}
%Also, near optimality for average cost policies as well...due to ACOE and \cite[Theorem 7.3.4]Lecture Notes. }

\section{Existence Results for Optimal Policies: Average Cost Criterion}

\subsection{Existence under strong sense admissible policies}\label{SconditionAvg}

%\sy{Implications for Average Cost: REFINE: Using the above for average cost requires the following: either minorization, which is not possible, or Wasserstein regularity of the kernel and the $\bar{c}$: \[E\bigg[c(X_k,u) \bigg| (Y_{(-\infty,k]}, U_{(-\infty,k-1]}=s) \bigg] =: \bar{c}(s,u),\]
%is to be Lipschitz in the state. If we have this, and that the kernel $\mathbb{P}^d$ is Wasserstein continuous, then the above also applies for the average cost problem. EXPAND THIS with references...}

\sy{
As a prelude result, we first note that we can have the following counterpart to Theorem \ref{LipschitzValueFInf} for the average cost problem. 
\begin{theorem}\label{avgCostLips}
\begin{itemize}
\item[(i)] Suppose that Assumptions \ref{contCond}, \ref{contConditinal} and \ref{Wasser1} hold with $\bar{K}_2 < 1$. Then, there exists a solution to the average cost optimality equation.
\item[(ii)] Suppose that Assumptions \ref{contCond} and \ref{contConditinal} hold. Then, there exists an optimal invariant measure and an associated optimal control policy\footnote{As in Theorem \ref{Belief-MDPMethod} , the optimality result may only hold for a restrictive class of initial conditions or initializations.}.
\end{itemize}
\end{theorem}

\textbf{Proof.} 
(i) See \cite[Lemma 2.2]{demirci2023geometric}, via the vanishing discount method. (ii) This result follows the proof of Theorem \ref{Belief-MDPMethod} given the weak continuity under Assumptions \ref{contCond} and \ref{contConditinal}.   
\qed

In the following, we present more relaxed conditions, tailored to the average cost criteria.

}
\subsection{Existence under conditionally state independent control policies}\label{RelaxedWSA}

%Once we have established controlled Markov models, we now discuss the classes of control policies considered under the %infinite dimensional controlled Markov model presented. 

We recall that in the theory of stochastic control, to facilitate stochastic analysis for arriving at existence, structural or approximation results (e.g. via continuity-compactness properties), the set of control policies may be enlarged. This is often referred to as a {\it relaxation} of control policies. Relaxations have been very effective in optimal control, with a very prominent example being Young measures in deterministic optimal control \cite{young1937generalized}. These allow one to use weak topologies on the spaces of probability measures to study existence, optimality, and structural results, rather than working with the space of measurable functions only (whose compactness conditions under an appropriate metric would be too restrictive). A key aspect of such relaxations is that any relaxation should not allow for optimal expected cost values to be improved; they should only be means to facilitate stochastic analysis. We can classify some policies and various relaxations as follows:

%\begin{enumerate}
%%\item Policies which depend only on the past information $I_t$ at any time $t$, i.e., measurable on the sigma field generated by the local information, are known as {\it strict-sense} admissible policies as in (\ref{eq_control}). {\it Randomized admissible} policies are those that are admissible based on the local history, but the control action set is the set of probability measures. Equivalently \cite[Lemma 1.2]{gihman2012controlled} \cite[Lemma 3.1]{BorkarRealization}, there is independent randomization at each time stage: For every $t \in \mathbb{Z}_+$, $u_t = \gamma_t(I_t,R_t)$ where $R_t$ is an independent $[0,1]$-valued random variable. %In particular, $(\ref{eq_control})$ is replaced with
%%%\begin{equation} \label{eq_control2}
%%%\gamma_t: I_t \mapsto u_t \in {\cal P}(\mathbb{U}),\quad t\in \mathbb{Z}_+
%%%\end{equation}
%%\item As noted earlier, policies that are based on the belief state, under a belief-MDP reduction, are known as {\it separated} policies. As noted, such policies are known to be optimal for a wide class of problems whenever measurable selection criteria (e.g. \cite[Chapter 3]{HernandezLermaMCP}) can be established. %We note that for linear quadratic models, separation often has a more specialized meaning, where one considers only those policies which depend on the information process $I_t$ only through the conditional expectation of the state $E[X_t| I_t]$  \cite{ber00,KushnerKalmanFilter,Caines}. 
%\item
(i) {\it Wide sense admissible policies} introduced by Fleming and Pardoux \cite{FlPa82} and prominently used to establish the existence of optimal solutions for partially observed stochastic control problems. Borkar \cite{Bor00,Bor07,Bor03} (see also Borkar and Budhiraja \cite{BoBu04}) have utilized these policies for a coupling/simulation method to arrive at optimality results for average cost partially observed stochastic control problems. The approach here is to first apply a Girsanov type (see Borkar \cite{Bor00,Bor07} for discrete-time models and Witsenhausen \cite{wit88} for decentralized stochastic control where the change of measure argument leads to {\it static reduction}) transformation to decouple the measurements from the system via an absolute continuity condition of measurement variables conditioned on the state, with respect to some reference measure; and then use independence properties: In the discrete-time case, $\{Y_n\}$ is i.i.d. and independent of $X_0$ and the system noise $\{W_n\}$, and $\{U_0,\ldots,U_n,Y_0,\ldots,Y_n\}$ is independent of $\{W_n\}$, $X_0$, and $\{Y_m, m>n\}$, for all $n$. 
%As cautiously noted in \cite[Section 8]{saldiyukselGeoInfoStructure}, if we only required that $\{U_n,Y_0,\ldots,Y_n\}$ is %independent of $\{W_n\}$, $X_0$, and $\{Y_m, m>n\}$, for all $n$, then a counterexample on the validity of the %relaxation can be established building on results from the quantum information theory literature  \cite{CHSH1969}.

(ii) Policies defined by conditional independence (as in the dynamic programming formulation \cite{YukselWitsenStandardArXiv} for decentralized stochastic control). In the context of our setup here, the actions are those that satisfy $U_n$ being conditionally independent of the past $(X_0,W_{[0,n]},V_{[0,n]})$ given the local information $\{Y_{[0,n]}, U_{[0,n-1]}\}$:
\begin{align}\label{CHIP}
U_n \leftrightarrow \{Y_{[0,n]}, U_{[0,n-1]}\} \leftrightarrow \{X_0,W_{[0,\cdots)},V_{[0,\cdots)}\}
\end{align}
By the notation $A \leftrightarrow B \leftrightarrow C$ for three random variables $A, B, C$ defined on a common probability space, we mean that $A$ and $C$ are conditionally independent given $B$; that is, $\mathrm{Prob}(A \in \cdot | B, C) = \mathrm{Prob}(A \in \cdot | B)$ almost surely.

We note that, here an absolute continuity condition required for a Girsanov-type measure transformation is not necessary apriori. We will call such policies {\it Conditionally-Exogenous Variable-Independent Policies} since the conditional independence holds between action, information, and the exogenous random variables in the system. 

If we had independent static reduction, with measurements $\tilde{Y}_n$, by considering the method for static measurements with expanding information structure (this is the information structure that one would obtain for a POMDP under the absolute continuity conditions) in \cite[Theorem 5.6]{YukselSaldiSICON17} or \cite[Theorem 4.7]{saldiyukselGeoInfoStructure} by expressing all the cost-relevant uncertainty in terms of $\tilde{Y}_n$, we have the following version of (\ref{CHIP})
\begin{align}\label{LMRUDP2}
U_n \leftrightarrow \{\tilde{Y}_{[0,n]}, U_{[0,n-1]}\} \leftrightarrow \{\tilde{Y}_{[n+1,\cdots)}\}
\end{align}

%Note that this is precisely the condition given in \cite{Bor00,Bor07,Bor03} (see also \cite{BoBu04}). 

%Thus, we can conclude that the Conditionally-Exogenous Variable-Independent Policies generalize Wide-Sense Admissible Policies in that they are equivalent when the Girsanov/Borkar/Witsenhausen reduction applies, but do not necessarily require such a reduction. The above also shows that such policies strictly generalize admissible randomized policies. \sy{Is this similar to how one can apply such wide-sense-admissible policies for degenerate diffusions? There one can still define them, without absolute continuity. Borkar's 2005 survey paper touches on this, but read this a bit more. There, policies are assumed to be measurable already with respect to the Brownian motion somehow, so open-loop in the first place...Setup is different.}
%% The policies can only use the assumed policies. 
%%Thus, different from conditionally-exogenous variable-independfent Policies, first the system is reduced to a static form with independent measurements, and the information available is assumed to 
%%the last item above can be considered to be those which first admit a static reduction/independent-measurements reduction, and then policies are measurable on the new process. 
%\end{enumerate}

In our paper, we will consider a further relaxation (to be called {\it Conditionally-State-Independent Policies}) which satisfies instead of (\ref{CHIP})
\begin{align}\label{CSIP}
U_n \leftrightarrow \{Y_{[0,n]}, U_{[0,n-1]}\} \leftrightarrow X_n,
\end{align}
together with (\ref{policyIndepEx}) below, to facilitate the infinite horizon analysis. \sy{Note that, the above does not imply conditional independence from past state realizations; however, in our analysis we will find that this relaxation is sufficient to arrive at an optimal policy; that is, this is an attainable relaxation.}

If the process continues on since the indefinite past, we write (\ref{CSIP}) with
\begin{align}\label{CSIPP}
U_n \leftrightarrow \{Y_{(-\infty,n]}, U_{(-\infty,n-1]}\} \leftrightarrow X_n
\end{align}

The above, however, is also to satisfy the condition that the process $(Z_k,U_k)$ is a controlled Markov chain as given in Lemma \ref{CMCZ_k} such that for all policies $\gamma$
\begin{align}\label{policyIndepEx}
P^{\gamma}(Z_{k+1} \in \cdot | Z_k=z_k, U_k=u_k) &:= P(Z_{k+1} \in \cdot | Z_k=z_k, \gamma(z_{[0,k]})=u_k) \nonumber \\
&= \mathbb{P}^a(\cdot | Z_t=z, U_t=u) 
\end{align}
This ensures that {\it policy-independence of conditional expectation} \cite{witsenhausen1975policy} holds, which is a required condition for the convex analytic method (as well as classical dynamic programming)\footnote{In the absence of this natural condition, one can construct counterexamples: Let $x_{t+1}=x_t + u_t + w_t$ where $w_t$ is i.i.d., and $y_t=x_t$; with the stage-wise cost $x_t^2$ . Then, under the policy $u_t = -w_t$, a better cost than a strictly admissible one can be attained; this policy would satisfy (\ref{CSIPP}) but not (\ref{policyIndepEx}).}. 

Thus, we require that (\ref{CSIPP})-(\ref{policyIndepEx}) hold for {\it Conditionally-State-Independent Policies}).

%While \cite{FlPa82} has arrived at existence results for POMDPs and \cite{Bor00,Bor07,Bor03} have arrived at optimality results for average cost stochastic control via separated designs, non-separated policies have also been studied: for continuous-time models, Bismut \cite{bismut1982partially} has arrived at  existence results, through an approach which avoids separation. As reviewed earlier, in discrete-time, \cite[Section 5.4.2, Theorem 5.6]{YukselWitsenStandardArXiv} presents such a characterization, without reducing the problem to a belief-MDP.
%
%\begin{remark}
%Under the conditions to be stated in the paper, conditionally-state-independent relaxation is an acceptable relaxation (in that the optimal performance does not improve): We can view the control policy as a {\it glueing} of $(I_t,X_t)$ to $(I_t,U_t)$ to form a Markov chain $U_t \leftrightarrow I_t \leftrightarrow X_t $. $(I_t,U_t)$ admits a convex representation (see Lemma 1.2 in \cite{gihman2012controlled}, or Lemma 3.1 of \cite{BorkarRealization}) so that $U_t = \gamma_t(I_t,R_t)$. The question now is whether $R_t$ is independent of all exogenous random variables: to answer this question, we will first obtain a lower bound by only imposing the Markov chain condition (\ref{CSIP}). Once we obtain an optimal solution (as a measure), we will show that without any loss, we can take $R_t$ to be independent from $X_t$ to realize the optimal measure. %Thus this relaxation will reduce to having independent randomization and this is an acceptable relaxation. \qed
%\end{remark}

\subsection{On the existence of an optimal stationary policy}

We first present a supporting result in the following, which will be crucial in our analysis to follow.

%\subsubsection{Weak closedness of sets of probability measures satisfying conditional independence under stationarity or invariant conditional %probabilities}

\begin{lemma}\label{keyLemmaClosed}
a) Let $Y^1_n,Y^2_n,U^2_n$ be a sequence of random variables that satisfies for every $n$, the conditional independence property: 
\[Y^1_n \leftrightarrow Y^2_n \leftrightarrow U^2_n,\]
with joint measure $P_n$, where the marginal on $(Y^1_n, Y^2_n)$ is fixed throughout the sequence. In this case, if $P_n \to P$ weakly, the limit measure $P$ also satisfies
\[Y^1 \leftrightarrow Y^2 \leftrightarrow U^2\]

b) Let $Y^1_n,Y^2_n,U^2_n$ be a sequence of random variables that satisfies for every $n$: 
\[Y^1_n \leftrightarrow Y^2_n \leftrightarrow U^2_n\]
with measure $P_n$ where the conditional probability $P_n(Y^1_n \in dy^1 | Y^2_n=y^2) = \kappa(dy^1|y^2)$ is fixed throughout the sequence and that $\kappa(dy^1|y^2)$ is a weakly continuous kernel (i.e. $y^2 \mapsto \int f(y^1) \kappa(dy^1|y^2)$ is continuous and bounded for every bounded continuous $f$). In this case, if $P_n \to P$ weakly, then, the limit measure also satisfies
\[Y^1 \leftrightarrow Y^2 \leftrightarrow U^2\]

\end{lemma}

\textbf{Proof.} See Section \ref{Proofkeylemmaclosed}. \qed

%We now present a slightly relaxed version which will be important for the controlled setup.
%
%\begin{lemma}\label{keyLemmaClosedCont}
%Let $Y^1_n,Y^2_n,U^2_n$ be a sequence of random variables that satisfies for every $n$: 
%\[Y^1_n \leftrightarrow Y^2_n \leftrightarrow U^2_n\]
%with measure $P_n$ where the conditional probability $P_n(Y^1_n \in dy^1 | Y^2_n=y^2) = \kappa(dy^1|y^2)$ is fixed throughout the sequence and that $\kappa(dy^1|y^2)$ is a weakly continuous kernel. In this case, if $P_n \to P$ weakly, then, the limit measure also satisfies
%\[Y^1 \leftrightarrow Y^2 \leftrightarrow U^2\]
%\end{lemma}
%
%\textbf{Proof.} See Section \ref{Proofkeylemmacontrolclosed}. \qed

%\subsection{Unique ergodicity of nonlinear filters}
%The same results in the control-free case will also apply here: Unique ergodicity of $I_t$ and filter stability implies the unique ergodicity of the filter process $\pi_t$. State as a theorem and prove it. 

We now present a key intermediate result.
\begin{lemma}\label{closednessTight2}
Let $\mathbb{X}, \mathbb{Y}, \mathbb{U}$ be compact. Under Assumptions \ref{contCond}, \ref{contConditinal} and \ref{ForgetCondition}, the set of invariant occupation measures
\begin{eqnarray}\label{ergodInv333}
{\cal G} &=& \{v \in {\cal P}(\mathbb{X} \times \mathbb{Y}^{\mathbb{Z}_-} \times \mathbb{U}^{\mathbb{Z}_-}): \nonumber \\
&& v(B \times \mathbb{U}) = \int_{z,u} \mathbb{P}^a(z_{t+1} \in B|z,u) v(dz,du) , \quad B \in {\cal B}(\mathbb{X} \times \mathbb{Y}^{\mathbb{Z}_-} \times \mathbb{U}^{\mathbb{Z}_-}) \} \nonumber \\
\end{eqnarray}
 that simultaneously satisfy the Markov chain 
\[X_0 \leftrightarrow (U_{(-\infty,-1]}),Y_{(-\infty,0]}) \leftrightarrow U_0,\]
or, equivalently, that belong to
\begin{eqnarray}\label{ergodInv444}
{\cal H} &=& \{v \in {\cal P}(\mathbb{X} \times \mathbb{Y}^{\mathbb{Z}_-} \times \mathbb{U}^{\mathbb{Z}_-}): \nonumber \\
&& \quad  v(dx_0,dy_{(-\infty,0]},du_{(-\infty,-1]},du_0) \nonumber \\
&& \quad \quad  = v(dx_0,dy_{(-\infty,0]},du_{(-\infty,-1]}) P(du_0 | dy_{(-\infty,0]},du_{(-\infty,-1]})\}, 
\end{eqnarray}
is a weakly compact set.
\end{lemma}

\textbf{Proof.} We note that under weak continuity of the transition kernel for $\mathbb{P}^a$, ${\cal G}$ is a closed set under weak convergence. By Lemma \ref{keyLemmaClosed}(b), as the kernel
\[P(dx_t | y_{(-\infty,t]},u_{(-\infty,t-1]})\]
is weakly continuous (by Assumption \ref{contConditinal}) and time-invariant (by Assumption \ref{ForgetCondition}), ${\cal H}$ is also closed. Since the state space is compact, the set of probability measures on $(\mathbb{X} \times \mathbb{Y}^{\mathbb{Z}_-} \times \mathbb{U})$ is tight. The result follows since a closed subset of a tight set is weakly compact. \qed

Building on the above, we are able to state the following:

\begin{theorem}\label{controlledPOMDPE}
Let $\mathbb{X}, \mathbb{Y}, \mathbb{U}$ be compact and Assumptions \ref{contCond}, \ref{ForgetCondition} and \ref{contConditinal} hold. \sy{Then, within the class of conditionally-state-independent-policies (\ref{CSIP}), an optimal control policy exists; this policy is stationary and can be written as}:
\[U_t = \gamma(Y_{(-\infty,t]}, U_{(-\infty,t-1]}, R_t), \]
where $R_t$ is an i.i.d. sequence. 
\end{theorem}

\textbf{Proof.} See Section \ref{controlledPOMDPEProofSection}. \qed

%\sy{Without any loss, the policy can be taken to be deterministic -work with Ari Araposthasis via the convex analytic method-. No..the set of strategic measures is not convex here!!!!!}

%\sy{for a discounted cost horizon version of this problem, we can use finite memory policies, but we would need uniform convergence properties. Finite window policies would quantize the memory. Furthermore {\bf the quantization error would be uniform}, since the topology is not on the probability measure space as we do with Ali, but on the state space. Since we have compactness if the measurement space is compact, the result will go through. ADD AS A SEPARATE SECTION!!!!! WEAK FELLERNESS IS EVEN VERY USEFUL!!!!..}

\begin{remark}
The above is a relaxation on the known results where one needs to reduce the POMDP to a belief-MDP which is weak Feller. As noted earlier, the weak Feller conditions require total variation continuity on either the system kernel or the measurement kernel. Only weak continuity of both the transition kernel and the measurement kernel is needed for the result above, though filter stability is also imposed.
\end{remark}

%\sy{1st benefit: So far, for this section says one new result. We have a far weaker condition; we only require the weak continuity conditions on the kernels and the measurement kernel. The conditions that require weak continuity of the belief-MDP are somewhat stronger, as far as the sufficient conditions are concerned!  We require compactness conditions, however, on the state space and action space. Nonetheless, this is a very important relaxation.}

\subsection{On the initial state distribution realization for an optimal stationary policy}

%\sy{Perhaps a separate discussion for discounted cost problems??}
In the above, there is an important operational question: The actual system is a one-sided process with $X_0 \sim \mu$ being the starting variable. To address this issue, in the following we explain that we can randomly generate the past and use the past realizations to apply optimal control. Under some conditions, the average cost will converge to the optimal cost. 

Our strategy on the initialization is as follows: (i) We will first show that an optimal invariant measure can be taken to be ergodic (i.e., an invariant measure which cannot be expressed as a convex combination of multiple invariant measures). (ii) Then we will construct an initialization which is in the attractor set of this optimal invariant measure, under a mild absolute continuity condition. We will generate the measurement and control actions $Y_{(-\infty,-1]},U_{(-\infty,-1]}$ according to a stationary measure, and then $X_0 \sim \mu$, $Y_0$ accordingly, and continue with the process according to the optimal policy under Theorem \ref{controlledPOMDPE}. 

\begin{lemma}\label{ergodicInv}
Without any loss of optimality, an optimal invariant occupation measure can be assumed to be ergodic.
\end{lemma}

\textbf{Proof.} \sy{
Let $\nu$ be an optimal occupation measure, which leads to a policy $\gamma$. Under $\gamma$, the {\it state} process 
\[(Y_{(-\infty,k]},U_{(-\infty,k-1]},X_k),\] 
is a Markov chain with the marginal of the invariant measure $\nu$ on this state being invariant. By an ergodic decomposition theorem, every stationary measure can be expressed as a convex combination of ergodic invariant measures with disjoint supports in the following sense: a stationary measure is a barycenter of the ergodic measures, i.e., if $\mu \in {\cal S}$ := the set of stationary measures and $E :=$ the set of its extreme points, i.e., the ergodic measures, then
\[ \int f(x) \mu(dx) = \int \int f(x)\nu(dx) \Psi(d\nu) \]
for every bounded measurable $f$ and a measure $\Psi$ on $E$. Thus a stationary
measure can be viewed as a random mixture of the ergodic measures with mixing probability $\Psi$. Thus, let $\nu$ be expressed as a convex combination of measures $\nu_{\beta}$, parametrized by $\beta$, where each of them are also {\it invariant} and {\it ergodic}. Each of the $\nu_{\beta}$ measures satisfies (\ref{ergodInv333}) and (\ref{ergodInv444}): The first holds by the invariance of each of these measures. The latter holds by definition: conditioned on the support of each $\nu_{\beta}$, $\nu$ satisfies (\ref{ergodInv444}) (note that the control policy is fixed) and hence this holds also for $\nu_{\beta}$.}

%Then, the infimum over all these ergodic invariant measures will be at least as good as the performance of $\langle \nu, c \rangle$, and each of the ergodic invariant measures $\nu_{\beta}$ will also simultaneously belong to both ${\cal H}$ and ${\cal G}$, by definition since the control policy is fixed for $\nu$, and will be the same for each of the invariant measure. That is for the constrained problem involving both (\ref{ergodInv333}) and (\ref{ergodInv444}), by the invariance property, the ergodic measures in any ergodic decomposition will also simultaneously satisfy both equations. [EXPLICIT ENOUGH? PERHAPS EXPRESS THE OPTIMIZATION PROBLEM MORE EXPLICITLY TO SHOW THAT THE INVARIANT ERGODIC MEASURES SATISFY BOTH CONDITIONS]
Since the cost in (\ref{expCostGH}):
\begin{eqnarray}
\inf_{v \in {\cal G} \cap {\cal H}} \int v(dx,dy_{(-\infty,0]},du_{(-\infty,-1]},du) c(x,u)
\end{eqnarray}
is linear in the space of (signed) measures, we conclude that without any loss we can take $\nu$ to be from the set of ergodic invariant measures.
% \sy{Does the intersection being a constraint lead to an issue? extremen points may not belong to ${\cal H}$? If does, then no problem...no, it is OK.}
\qed

The above then leads to our realization result, whose proof is in the appendix. 

\begin{theorem}\label{FilterInvarianceMeasure2}
Let Assumptions \ref{contConditinal} and \ref{ForgetCondition} hold. Let $U_t = \gamma(Y_{(-\infty,t]},U_{(-\infty,t-1]},R_t)$ be an optimal stationary policy with $R_t$ being an i.i.d. noise process. Let \[g(y_{(-\infty,k]},u_{(-\infty,k-1]},x_k):= E[c(x_k, \gamma(y_{(-\infty,k]},u_{(-\infty,k-1]},R_k))]\]
 where the expectation is over $R_k$. Let $\bar{P}$ denote an (invariant) process measure on \[(Y_{(-\infty,k-1]},U_{(-\infty,k-1]},Y_k, X_k),\] and let $P_0$ be the marginal of the invariant measure on $\{Y_{(-\infty,-1]},U_{(-\infty,-1]}\}$ and $\mu$ be the prior measure on $X_0$ as imposed (i.e., given in the problem statement) on the controller/decision maker. If $P_0 \times \mu \ll \bar{P}$ (note that the distribution on $X_0$ specifies that on $Y_0$); that is, absolute continuity holds when the (past history) initialization is independent of the state process initialization, then
\begin{eqnarray}
&& \lim_{T \to \infty} {1 \over T} E_{P_0 \times \mu} \bigg[\sum_{k=0}^{T-1} c(X_k, \gamma(Y_{(-\infty,k]},U_{(-\infty,k-1]},R_k)) \bigg] \nonumber \\
&& \quad \quad \quad \quad  = \int g(y_{(-\infty,0]},u_{(-\infty,-1]},x_0) \bar{Q}(dy_{(-\infty,0]},du_{(-\infty,-1]},x_0), \label{filterInitialization}
\end{eqnarray}
for some invariant $\bar{Q}$, which would necessarily be $\bar{P}$ under the assumed ergodicity of $\bar{P}$ via Lemma \ref{ergodicInv}, in which case the optimal control/decision cost will be attained.
%Thus, if we generate $\pi_{-1}$, the filter process at time $-1$ randomly according to some invariant measure for $\pi_t$ consistent with $\bar{P}$, then under the condition stated, almost surely, the filter will correct itself asymptotically.
\end{theorem}

%\begin{proposition}
%Even with an incorrect initialization $P_0$ on $Y_{(-\infty,-1]}$, the process $(Y_{(-\infty,t]},X_t)$ with initial measure $P_0 \times \pi_X$ is a Markov process with transition kernel $\mathbb{P}$.
%\end{proposition}

\begin{remark}
Note that, in the above, with $\bar{P}$ and $P_0$ both stationary, we have that $\bar{P}$ decomposes as $P_0(dy_{(-\infty,-1]},du_{(-\infty,-1]})  \bar{P}(dx_0 | y_{(-\infty,-1]},u_{(-\infty,-1]})$. 
A sufficient condition for the aforementioned absolute continuity condition then is, $\bar{P}$ a.e.,
\[\mu(dx_0) \ll  \bar{P}(dx_0 | dy_{(-\infty,-1]},du_{(-\infty,-1]})\]
\end{remark}

\begin{remark}
Observe that in the analysis above we constructed $y_{(-\infty,-1]}, u_{(-\infty,-1]}$ and not $y_{(-\infty,0]}$, since $\mu$ on $X_0$ induces directly a measure on $Y_0$. Note that with $P_0 \times \mu$ viewed as an initial measure, when the process evolves, the evolution is consistent with the actual realizations under the true measure. That is, {\it the initialization does not affect the evolution of the measure}. The reason is that the process $X_0$ determines the probability of the future events, independent of the history initialization via:
\[\mathbb{P}^a(dz_{t+1} | z_t,u_t) = Q(dy_{t+1}|x_{t+1}) P(dx_{t+1} | x_t,u_t) P(dy_{(-\infty,t]} | y_{(-\infty,t]} ) P(du_{(-\infty,t]} | u_{(-\infty,t]} ).\]
%is such that $Y_t$ is generated according to both the true process and the actual process. 
\end{remark}

\section{The Control-Free Case: Implications on Unique Ergodicity of Nonlinear Filters and Existence of an Optimal Stationary Policy}\label{conFreeCase}
In this section, we interpret some of the results presented in the controlled-case in the context of the control-free case. This serves both as an application and validation of the approach, but also presents new results on the nonlinear filtering problem. 

Even though the control-free case can be viewed as an instance of the controlled case, here some of the assumptions can be relaxed. Accordingly, we will state two results. One is on the stationarity properties on the optimal decision policies, and the other is on how to realize it by selecting the initial prior appropriately. 

We have the following assumption for some of the results to follow. This will be a complementary condition.
\begin{assumption}\label{statX}
$\{X_t\}$ is stationary.
\end{assumption}

Note that if $X_t$ is stationary, so is the pair process $(X_t, Y_t)$. By a standard argument (e.g. {\sl Chapter 7} in \cite{durrett2010probability}), we can embed the one-sided stationary process $\{X_k, k \in \mathbb{Z}_+\}$ into a bilateral (double-sided) stationary process $\{ X_{k}, k \in \mathbb{Z} \}$. 

%The following is the main result of this section

\begin{theorem}\label{MeanOptimal}
Let $\mathbb{X}, \mathbb{Y}, \mathbb{U}$ be compact and Assumption \ref{contCond} hold. Under either Assumption \ref{statX} or Assumptions \ref{contConditinal} and \ref{ForgetCondition} an optimal decision/control policy exists. This policy is stationary.
\end{theorem}

\textbf{Proof.} See Section \ref{SMeanOptimalProof}. \qed 

\begin{corollary}\label{FilterInvarianceMeasure222}
Let either (i) Assumption \ref{statX}, or, (ii) under the hypotheses in Theorem \ref{MeanOptimal}, Assumption \ref{assumptionMC}, hold. Then, the filter process admits an invariant probability measure.
\end{corollary}

% note that by Assumption \ref{assumptionMC}, we have a measurable function of the past...so invariance of the past implies that of the filter. note that there is
% no control here....

As noted above, in the uncontrolled setting, \cite{bhatt2000markov} and \cite{budhiraja2002invariant} have established weak continuity conditions (i.e., the weak-Feller property) of the nonlinear filter process (i.e., the belief process) in continuous time and discrete time, respectively; where the total variation continuity of the measurement channel was imposed. These results were used to establish the existence of an invariant probability measure for the belief process. The above shows that for the existence of an invariant probability measure, one may not need to invoke the continuity conditions.

%\begin{remark}
%Note that we could view:
%\[E[f(X_t) | Y_{(-\infty,t]}]\]
%as the solution to the optimization problem with cost
%\[c(X_t, U_t) = \|f(X_t) - U_t\|^2\]
%Under any continuous and bounded $f$, we can obtain a solution; this will be stationary and will lead to the existence of an invariant probability measure for the filter process, since a countable collection of such $f$ functions will determine the filter process. 
%% and it will lead to an ergodic process $X_t, Y_{(-\infty,t]}$; the solution will fully characterize the ergodic process $\pi_t$. This will ensure the uniqueness of an invariant probability measure, through some additional measure theoretic analysis (of representing continuous and bounded functions on the space of prob. measures). Stettner and Di Masi have some related results \cite{stettner1991invariant, di2008ergodicity}. It would be interesting to establish the equivalence of the ergodicity results obtained earlier in the paper; see Theorem \ref{uniqueInvMeasure}.
%\end{remark}

A corollary (to Theorem \ref{MeanOptimal} and Lemma \ref{measurabilityBelief}) is the following. 
\begin{corollary}\label{FilterInvarianceMeasure22}
Let $\mathbb{X}, \mathbb{Y}, \mathbb{U}$ be compact and Assumption \ref{contCond} hold. Additionally, let the hypotheses in Theorem \ref{MeanOptimal} hold. If $Y_t$ is uniquely ergodic and if Assumption \ref{assumptionMC} holds, then the filter is uniquely ergodic.
%Let Assumption \ref{statX} hold. Then, the filter process admits an invariant probability measure.
\end{corollary}

\textbf{Proof.} We have seen that, by the hypotheses in the statement of Theorem \ref{MeanOptimal}, there exists at least one invariant probability measure. Suppose that there were two distinct invariant probability measures, $\eta_1$ and $\eta_2$. Since a countable collection of continuous and bounded functions can be used to distinguish probability measures, we will consider $\int \eta_i(d\pi) f(\pi)$ for such a countable collection of functions $f$, for $i=1,2$. Consider the joint process $(\pi_t, Y_{(-\infty,t]})$ with invariant measure $\kappa_i$ and with a marginal invariant measure on $\pi_t$ as $\eta_i$. Since the marginal on $Y_{(-\infty,t]}$ is uniquely ergodic, for every invariant measure $\kappa_i$ on the joint process, the marginal will be a constant measure, call $\psi$. Therefore:
\[\int \eta_i(d\pi) f(\pi) =  \int \kappa_i(d\pi, dy_{(-\infty,t]}) f(\pi) = E_{\kappa_i}[f(\pi)] = E_{\kappa_i} [ E_{\kappa_i}[f(\pi) | Y_{(-\infty,t]}] ] \]
However, $E_{\kappa_i}[f(\pi) | Y_{(-\infty,t]}] = E_{\kappa_i}[g(Y_{(-\infty,t]})] = E_{\psi}[g(Y_{(-\infty,t]})] $ for some measurable $g$, by Assumption \ref{assumptionMC}. Therefore, as this argument applies for any $f$ from the distinguishing family, we must have that $\eta_1=\eta_2$. \qed
 
Note that for $Y_t$ to be uniquely ergodic, we don't require $X_t$ to be uniquely ergodic. On the other hand, if $X_t$ is uniquely ergodic, $Y_t$ must be uniquely ergodic as well. For sufficient conditions on unique ergodicity of infinite-memory processes such as $Y_t$, see \cite[Theorem 3.8]{yukselSICON2017}; in particular the existence of an accessible element and a continuity condition (instead of weak continuity, where setwise continuity is imposed in Assumption \ref{contConditinal}; e.g., \cite[Corollary 4.2]{le2004stability} implies this condition as well).

Notably, building on \cite[Theorem 2]{DiMasiStettner2005ergodicity} and \cite[Prop 2.1]{van2009uniformSPA}, it can be shown, via a functional analytic argument, that almost sure filter stability in the total variation sense (which can be relaxed to weak merging) and the uniqueness of an invariant probability measure on the state process leads to unique ergodicity of the filter process. It thus turns out that a complementary condition, and with an alternative argument (e.g. to \cite[Theorem 2]{DiMasiStettner2005ergodicity}), can also be established: unique ergodicity of the measurement process and a measurability condition related to filter stability leads to unique ergodicity.

\begin{remark}\label{RemarkContFree}  For the problem at hand, for the control-free setup, it is evident that the decision maker can always apply $\arg \min E[c(X_k,U_k) | Y_{[0,k]}]$ at any $k \in \mathbb{N}$. This policy will lower bound the expected cost under any policy, including the aforementioned policy. The conclusion is that, asymptotically, they are equivalent.
The main novelty here is the optimality of a {\it stationary} policy for an average cost problem. 
\end{remark}

{\bf On the realization problem for an optimal stationary policy.} Assumption \ref{statX} will not be applicable if the initialization of the process is to be arbitrary. Let $\bar{P}$ be the (unique) invariant probability measure on $\{Y_{(-\infty,t]}\}$. If the initial measure $P_0$ is such that $P_0 \ll \bar{P}$, then the ergodic theory of Markov chains (see Theorem \ref{convergenceT3}(ii)) implies that for every measurable bounded $g$:
\[ \lim_{T \to \infty} {1 \over T} E_{P_0}[\sum_{k=0}^{T-1} g(Y_{(-\infty,k]})] = \int g(y_{(-\infty,0]}) \bar{P}(dy_{(-\infty,0]}) \]

%%In fact, this convergence is uniform over such $g$. 
%%Since the filter output for a given event $B$ is a measurable and bounded function of the past and current measurements, this also implies that
%%\[ \lim_{T \to \infty} {1 \over T} E_{P_0} \sum_{k=0}^{T-1} \pi_k(y_{(-\infty,k]})(A) = \int \pi_k(y_{(-\infty,0]})(A)  \bar{P}(dy_{(-\infty,0]}) \]
%%with this convergence also being uniform over such $A$ as well as $k$. This uniformity may be crucial for us.
%%
%%\sy{We don't know if the filter is a measurable outcome of the past measurements. This in fact must be stated as an assumption; the initialization at the distant past may still be there. This must be stated as an assumption; unless the initial filter is assumed to be given.}
%%
%%%Another question then is how to generate $\pi_0$? }

We state the following, building on Theorem \ref{convergenceT3}(ii), and the proof method of Theorem \ref{FilterInvarianceMeasure2}.

\begin{theorem}\label{FilterInvarianceMeasure12}
Let Assumptions \ref{ForgetCondition} and \ref{contConditinal} hold. Let $u_t = \gamma(y_{(-\infty,t]},R_t)$ be an optimal stationary policy (by Theorem \ref{MeanOptimal}) with $R_t$ being an i.i.d. noise process. Let $g(y_{(-\infty,k]},x_k):= E[c(X_k, \gamma(Y_{(-\infty,k]},R_k))]$ where the expectation is over $R_k$. Let $\bar{P}$ denote an invariant process measure on $(Y_{(-\infty,k-1]},X_k)$ and let $P_0$ be the invariant measure on $\{Y_{(-\infty,-1]}\}$ and $\mu$ be the prior measure on $X_0$ as imposed by the controller/decision maker. If $P_0 \times \mu \ll \bar{P}$; that is, if the (filter) initialization is independent of the state process initialization and under the absolute continuity condition,
\begin{eqnarray}
 \lim_{T \to \infty} {1 \over T} E_{P_0 \times \mu}[\sum_{k=0}^{T-1} g(Y_{(-\infty,k]},X_k)] = \int g(y_{(-\infty,0]},x_0) \bar{Q}(dy_{(-\infty,0]},x_0). \label{filterInitialization0}
\end{eqnarray}
%with the convergence being uniform over $g$ and thus for the filter events. 
The optimal control/decision cost will be attained.
%Thus, if we generate $\pi_{-1}$, the filter process at time $-1$ randomly according to some invariant measure for $\pi_t$ consistent with $\bar{P}$, then under the condition stated, almost surely, the filter will correct itself asymptotically.
\end{theorem}

%\sy{I noted here that: the filter process is not a measurable function of $y_{(-\infty,k]}$. Consider Ramon's example: $X_{k+1}=X_k \oplus 1$, and $Y_k$ an i.i.d. process.}
%
Observe again, as in the controlled setup, that in the analysis above we constructed $y_{(-\infty,-1]}$ and not $y_{(-\infty,0]}$, since $\mu$ on $X_0$ induces directly a measure on $Y_0$. Note that with $P_0 \times \mu$ viewed as an initial measure, when the process evolves, the evolution is the correct one consistent with the actual realizations under the true measure and thus the initialization does not affect the evolution of the measure: $X_0$ determines the probability of the future events: the transition kernel:
\[\mathbb{P}^a(dz_{t+1} | z_t) = Q(dy_{t+1}|x_{t+1}) P(dx_{t+1} | x_t) \delta_{y_{(-\infty,t]} }(dy_{(-\infty,t]})\]
is such that $Y_t$ is generated according to both the true process and the actual process. 

%That is,
%\[ P_{\bar{P}} (dz_{t+1} | z_t) =  P_{P_0 \times \mu} (dz_{t+1} | z_t) \]
%as long as the measure on $X_0$ is consistent with the actual system. 

%We state this formally as follows.
%
%\begin{proposition}
%Even with an incorrect initialization $P_0$ on $Y_{(-\infty,-1]}$, the process $(Y_{(-\infty,t]},X_t)$ with initial measure $P_0 \times \mu$ is a Markov process with transition kernel $\mathbb{P}$.
%\end{proposition}
%Does there always exist $\pi_0$ which satisfies $P_0 \times \pi_{0} \ll \bar{P}$? 

Note that if $\bar{P}$ and $P_0$ are both stationary, we have that $\bar{P}$ decomposes as $P_0(dy) Q(dx_0|y_{(-\infty,-1]})$. Thus, as in the controlled-case, what we need is that $\mu(dx_0) \ll Q(dx_0|y_{(-\infty,-1]})$ $P_0$ a.e., so that convergence in the sense of (\ref{filterInitialization0}) holds.

\section{Conclusion}
We studied the optimal control of partially observed Markov Decision Processes (POMDPs) without the common paradigm of belief-separation. We first revisited this approach and presented conditions for the existence of optimal policies.

We then defined a Markov chain taking values in an infinite dimensional product space with the history process serving as the controlled state process and a further refinement in which the control actions and the state process are causally conditionally independent given the measurement/information process. 
We provided new sufficient conditions for the existence of optimal control policies under the discounted cost and average cost infinite horizon criteria. In particular, while in the belief-MDP reduction of POMDPs, weak Feller condition requirement imposes total variation continuity on either the system kernel or the measurement kernel, with the approach of this paper only weak continuity of both the transition kernel and the measurement kernel is needed (and total variation continuity is not) together with regularity conditions related to filter stability. For the discounted cost setup, we also establish near optimality of finite window policies via a direct argument involving near optimality of quantized approximations for MDPs under weak Feller continuity, where finite truncations of memory can be viewed as quantizations of infinite memory with a uniform diameter in each finite window restriction under the product metric. For the average cost setup, we provided new existence conditions and also a general approach on how to initialize the randomness to establish convergence to optimal cost.  In the control-free case, our analysis lead to new and weak conditions for the existence and uniqueness of invariant probability measures for nonlinear filter processes, where we showed that unique ergodicity of the measurement process and a measurability condition related to filter stability leads to unique ergodicity.

%We then defined a Markov chain taking values in an infinite dimensional product space with control actions and the state process causally conditionally independent given the measurement/information process, under which we provided new sufficient conditions for the existence of optimal control policies. In particular, while in the belief-MDP reduction of POMDPs, weak Feller condition requirement imposes total variation continuity on either the system kernel or the measurement kernel, with the approach of this paper only weak continuity of both the transition kernel and the measurement kernel is needed (and total variation continuity is not) together with regularity conditions related to filter stability. For the average cost setup, we provided a general approach on how to initialize the randomness which we show to establish convergence to optimal cost. For the discounted cost setup, we established near optimality of finite window policies via a direct argument involving near optimality of quantized approximations for MDPs under weak Feller continuity, where finite truncations of memory are viewed as quantizations of infinite memory with a uniform diameter in each finite window restriction under the product metric. In the control-free case, our analysis entails new conditions for the existence and uniqueness of invariant probability measures for nonlinear filter processes, where we showed that unique ergodicity of the measurement process and a measurability condition related to filter stability leads to unique ergodicity.

\section{Acknowledgements}
We are thankful to Eugene Feinberg, Naci Saldi, Ali D. Kara and Y. Emre Demirci for their comments. 

\appendix
\section{Discussion and Sufficient Conditions on Assumptions}\label{OnAssumptions}
As a relatively extreme case, a fully observed model, viewed as a POMDP via the measurement kernel (\ref{MDPasPOMDP}) satisfies Assumptions \ref{contConditinal}, \ref{assumptionMC} and \ref{ForgetCondition}; as noted earlier, such a model cannot be studied within the existing results on POMDP theory. 

For the case with finite measurement and actions, Assumption \ref{contConditinal} is satisfied by almost sure filter stability, though with a uniformity condition over priors, in the following sense
\begin{align}\label{UniformWeakConv}
\sup_{\nu} \|\pi_{n}^{\mu}-\pi_{n}^{\nu}\|_{BL} \to 0~P^{\mu} a.s.,
\end{align} 
where $BL$ denotes the bounded Lipschitz norm (or any other weak convergence inducing metric can be used), and $\pi_{n}^{\mu}, \pi_{n}^{\nu}$ denote the filter realizations with initializations at $\pi_0=\mu$ and $\pi_0=\nu$, respectively. This condition (including uniformity) holds under a contraction analysis via the {\it Hilbert metric}, as shown in \cite[Corollary 4.2]{le2004stability}. This also applies to the controlled case as shown in \cite[Lemma 5.11]{demirciRefined2023}.

Furthermore, by modifying the proof of \cite[Lemma 3.5]{MYDobrushin2020} to place $\sup_{\nu}$ inside the expectation, \cite[Corollary 3.7]{MYDobrushin2020} also leads to this condition in view of Borel-Cantelli Lemma as noted in \cite[Remark 3.10]{MYDobrushin2020}.

The sufficiency of (\ref{UniformWeakConv}) builds on the following: Observe that for a continuous and bounded $f$
\begin{eqnarray}
&& E[f(X_0)|y^n_{(-\infty,0]},u^n_{(-\infty,-1]}] - E[f(X_0)|y_{(-\infty,0]},u_{(-\infty,-1]}]  \nonumber \\
&& = E[f(X_0)|y^n_{(-\infty,0]},u^n_{(-\infty,-1]}] - E[f(X_0)|(y,u)^n_{(-\infty,-N]},y_{(-N+1,0]},u_{(-N+1,-1]}] \nonumber \\
 \label{cond1} \\
&& + E[f(X_0)|(y,u)^n_{(-\infty,-N]},y_{(-N+1,0]},u_{(-N+1,-1]}]  - E[f(X_0)|y_{(-\infty,0]},u_{(-\infty,-1]}] \nonumber \\ 
\label{cond2} 
\end{eqnarray}

The uniform filter stability condition (\ref{UniformWeakConv}) will allow us to truncate the past finite window: For every $\epsilon >0$ and $\pi=P(X_n \in \cdot | y_{(\infty,n-},u_{(\infty,n-1})$ select $N$ so that, the effect of the history is uniformly, over all priors $P(X_{-N} \in \cdot | (y,u)^n_{(-\infty,-N]})$, is less than $\frac{\epsilon}{2}$, so that (\ref{cond2}) is bounded uniformly over all sequences prior to $-N$. 

We then apply a continuity argument for the first term, (\ref{cond1}):

%Now, for the second item (\ref{cond2}), one could establish for a fixed $N$ uniform convergence (over measurement signals $(y,u)^n_{(-\infty,-N]}$) of
%\[E[f(X_0)|y^n_{(-\infty,0]},u^n_{(-\infty,-1]}] - E[f(X_0)|(y,u)^n_{(-\infty,-N]},y_{(-N+1,0]},u_{(-N+1,-1]}] \]
%to zero. Note that this is implied with the uniform convergence
%\begin{eqnarray}
%&&\lim_{y^n_{(-N+1,0]},u^n_{(-N+1,-1]} \to y_{(-N+1,0]},u_{(-N+1,-1]}} \sup_{\pi_{-(N-1)}} \nonumber \\
%&& \bigg| E^{\pi_{-(N-1)}}[f(X_0) | y^n_{(-N+1,0]},u_{(-N+1,-1]}]  - E^{\pi_{-(N-1)}}[f(X_0) | y_{(-N+1,0]},u_{(-N+1,-1]}] \bigg|. \nonumber 
%\end{eqnarray}
%To see why this holds as a consequence of filter stability, consider an auxiliary probability measure $\kappa = \frac{1}{2} P(X_{-N} \in \cdot | (y,u)^n_{(-\infty,-N]}) + \frac{1}{2} P(X_{-N} \in \cdot | (y,u)_{(-\infty,-N]})$. Then, $\kappa$ almost surely we have that $P^{P(X_{-N} \in \cdot | (y,u)^n_{(-\infty,-N]})}(X_N \in \cdot | y_{(-N+1,0]},u_{(-N+1,-1]}) - P^{\kappa}(X_N \in \cdot | y_{(-N+1,0]},u_{(-N+1,-1]})$ to 0. The same applies for $P^{P(X_{-N} \in \cdot | (y,u)_{(-\infty,-N]})}(X_N \in \cdot | y_{(-N+1,0]},u_{(-N+1,-1]}) - P^{\kappa}(X_N \in \cdot | y_{(-N+1,0]},u_{(-N+1,-1]})$ converging to zero $\kappa$ almost surely. 

For this first term, in case the the measurements and actions are finitely valued, the desired result follows since for sufficiently large $n$, the first $N$ coordinates of measurement and actions $y_{(-N+1,0]},u_{(-N+1,-1]}$, will need to match to satisfy proximity under the product metric so that for all sufficiently large $n$: 
\[y^n_{(-N+1,0]},u^n_{(-N+1,-1]}=y_{(-N+1,0]},u_{(-N+1,-1]},\]
making (\ref{cond1}) zero.

For the case with continuous measurements and actions, conditions in \cite[Lemma 4.6]{MYRobustControlledFS} suffice, together with the uniform filter stability condition presented above. 

It can be shown that a {\it uniform} filter stability condition, such as \cite[Corollary 4.2]{le2004stability} (via a Hilbert metric approach) as well as \cite[Corollary 3.7]{MYDobrushin2020} (via Dobrushin's coefficient method), implies Assumption \ref{contConditinal} for the case where the measurement and action spaces are finite. 

%Our next assumption is the following. 
%\begin{assumption}\label{ForgetCondition}[A Stationarity / Stability Condition]
%For every Borel $B$, and with $y = (\cdots, y_{-2}, y_{-1}, y_0)$, $u = (\cdots, u_{-2}, u_{-1}, u_0)$, under any control policy, for almost every $(y, u)$, we have that
%\[P(X_t \in B| Y_{(-\infty,t]}=y, U_{(-\infty,t-1]}=u) = P(X_{t+1} \in B| Y_{(-\infty,t+1]}=y, U_{(-\infty,t]} = u)\]
%\end{assumption}

Assumption \ref{ForgetCondition} is also implied by an almost sure filter stability condition, e.g. \cite[Corollary 4.2]{le2004stability}. Suppose that $P^{\gamma}\bigg(X_t \in \cdot \bigg| \bigcap_{n=1}^{\infty} \sigma(Y_{( -\infty,t]},U_{( -\infty,t-1]})  \vee \sigma(\pi_{(-\infty,-n]} ) \bigg)$ is $\sigma(I_t)$-measurable in the sense that $P$ a.s. 
\begin{eqnarray}\label{stabilityFC}
&&P^{\gamma}\bigg(X_t \in \cdot \bigg| \bigcap_{n=1}^{\infty} \sigma(Y_{( -\infty,t]},U_{( -\infty,t-1]})  \vee \sigma(\pi_{(-\infty,-n]} ) \bigg)  \nonumber \\
&& \qquad \qquad \qquad \qquad \qquad \qquad =  P^{\gamma}(X_t \in \cdot | \sigma(Y_{( -\infty,t]},U_{( -\infty,t-1]} ) ) \nonumber \\
&& \qquad \qquad \qquad \qquad \qquad \qquad =  P(X_t \in \cdot | \sigma(Y_{( -\infty,t]},U_{( -\infty,t-1]} ) ), 
\end{eqnarray}
where the last equality implies policy independence of the conditional probabilities.

\sy{We note that the above is not a trivial condition, since the measure $P$ may induce multiple ergodic invariant measures and therefore the conditioning on a particular $\pi_{-n}$ realization may lead to different realizations for $\pi_t$.}

% SY: Here, we assume that the two-sided process measure $P$ is well-defined. In the applications, we will require that it is 
% the result of a stationary process (satisfying the discounted optimal or the average optimal problen). In this case it will be well-%defined. It will be a process generated by some stationary
% measure on the control and dynamics. In general, however, it is not clear how one would define a two-sided process.

As before, with $\pi_t$ a conditional probability on the state given a (distant) prior at time $-n \in \mathbb{Z}_-$, and the information since then (or, equivalently, the information even prior to $-n$), we have for all Borel $A \subset \mathbb{X}$
\[ \pi_t(A) = E[1_{X_t \in A} | \sigma(Y_{( -\infty,t]},U_{( -\infty,t-1]})  \vee \sigma(\pi_{(-\infty,-n]}] \]
This conditional probability is well-known to be policy independent.
%In view of the above, we state the following assumption.

%\sy{Almost sure filter merging, in weak convergence/merging, under any prior measure on the measurement process satisfies the above. Here is the reason. $\pi_t$ is completely characterized by testing it on a countable set of continuous functions. Therefore, almost sure convergence for a countable set of continuous test functions will suffice. But the above precisely gives this characterization. Thus, almost sure weak merging is a sufficient condition. NO NEED FOR TOTAL VARIATION. THE CONDITION HERE IS DIFFERENT FROM KUNITA/VAN HANDEL.} 

\begin{lemma}\label{measurabilityBelief}
Condition (\ref{stabilityFC}) implies Assumption \ref{assumptionMC}.
\end{lemma}

\textbf{Proof.} Observe that $\pi_t(A) = E[1_{X_t \in A} | \sigma(Y_{( -\infty,t]},U_{( -\infty,t-1]})  \vee \sigma(\pi_{(-\infty,-n]}]$ for every $n > 0$ and thus \[\pi_t(A) = \lim_{n \to \infty} E[1_{X_t \in A} | \sigma(Y_{( -\infty,t]},U_{( -\infty,t-1]})  \vee \sigma(\pi_{(-\infty,-n]}].\]
Now, $E[1_{X_t \in A} | \sigma(Y_{( -\infty,t]},U_{( -\infty,t-1]})  \vee \sigma(\pi_{(-\infty,-n]}]$ is a bounded backwards martingale sequence with respect to the decreasing filtration $\sigma(Y_{( -\infty,t]},U_{( -\infty,t-1]})  \vee \sigma(\pi_{(-\infty,-n]})$, $n \in \mathbb{N}$ and by the backwards martingale theorem \cite{durrett2010probability}, the above limit will converge to
\[P\bigg(X_t \in A \bigg| \bigcap_{n=1}^{\infty} \sigma(Y_{( -\infty,t]},U_{( -\infty,t-1]})  \vee \sigma(\pi_{(-\infty,-n]} ) \bigg)\]

This then implies, via (\ref{stabilityFC}), $\pi_t(A)$ is $\sigma(I_t)$-measurable for any given $A$.

%But each member of this sequence is measurable on $\bigcap_{n=1}^{\infty} \sigma(Y_{( -\infty,t]},U_{( -\infty,t-1]})  \vee \sigma(\pi_{(-\infty,-n]})$, since $\bigcap_{n=1}^{\infty} \sigma(Y_{( -\infty,t]},U_{( -\infty,t-1]}) \vee \sigma(\pi_{(-\infty,-n]})$ determines the filter $E[1_{X_t \in A} | \sigma(Y_{( -\infty,t]},U_{( -\infty,t-1]})  \vee \sigma(\pi_{(-\infty,-n]}]$ (as every $n$ does so). Therefore, so is the limit $\pi_t(A)$.

Recall the following which builds on Theorem 2.1 of Dubins and Freedman \cite{DubinsFreedman} and Proposition 7.25 in Bertsekas and Shreve \cite{BertsekasShreve}: 
%\begin{theorem}\label{DubinsFreedmanTheorem}
{\it Let $\mathbb{S}$ be a Polish space, ${\cal P}(\mathbb{S})$ be the set of probability measures under the weak convergence topology and $(M,{\cal M})$ be a measurable space. A function $F: (M,{\cal M}) \to {\cal P}(\mathbb{S})$ is measurable on ${\cal M}$ if for all $B \in {\cal B}(\mathbb{S})$
$(F(\cdot))(B): M \to \mathbb{R}$ is measurable on ${\cal M}$, that is for every $B \in {\cal B}(\mathbb{S})$, $(F(\pi))(B)$ is a measurable function when viewed as a function from $M$ to $\mathbb{R}$.} See also \cite[Proposition 7.26]{BertsekasShreve}. By this result, it follows that $\pi_t$ is measurable on $\bigcap_{n=1}^{\infty} \sigma(Y_{( -\infty,t]},U_{( -\infty,t-1]})  \vee \sigma(\pi_{(-\infty,-n]})$. But by condition (\ref{stabilityFC}), it is also $\sigma(I_t)$ measurable, and as the considered spaces are standard Borel, a functional representation via the measurable function $g$ follows.
\qed

%Further analysis and sufficient conditions for Assumptions \ref{contConditinal}, \ref{ForgetCondition} and \ref{assumptionMC} are presented in Section \ref{sufficientCondKey} below. 

We again note that (\ref{stabilityFC}) is essentially a filter stability condition: Indeed, the assumption above, in the control-free setup, is related to the statement in \cite[Theorem 3.1(2)]{chigansky2010complete}, building on \cite{kunita1971asymptotic}.

 In the control-free case; Assumption \ref{ForgetCondition} is implied by stationarity or both Assumptions \ref{assumptionMC} and \ref{ForgetCondition} hold under almost sure filter stability, see e.g. \cite{chigansky2009intrinsic,chigansky2010complete}. 
 
 Assumption \ref{assumptionMC} also holds via the Hilbert metric condition given in \cite[Corollary 4.2]{le2004stability}. Complementing this, for both controlled and control-free setups, conditions in \cite{MYDobrushin2020} (via \cite[Theorem 2, Part 2]{kleptsyna2008discrete} leading also to almost sure stability under total variation) based on Dobrushin's coefficients of the measurement channel and the controlled transition kernel leads to almost sure filter stability and accordingly Assumption \ref{ForgetCondition}. 

Let $\mathbb{X}, \mathbb{Y}, \mathbb{U}$ be finite, and each be endowed with the discrete metric. Consider $s^n \to s$ in the product topology, where in this case we slightly revise the metric (since the spaces are finite) to be as follows: With $s^n = (y^n_{(-\infty,0]},u^n_{(-\infty,-1]}), s = (y_{(-\infty,0]},u_{(-\infty,-1]})$, let
\[\bar{d}(s^n,s) = \sum_{k=0}^{\infty} \rho^{-k} d_{\mathbb{Y}}(y^n_{-k},y_{-k}) + \sum_{k=1}^{\infty} \rho^{-k} d_{\mathbb{U}}(u^n_{-k},u_{-k})\]
for some  $\rho > 1$. 
 This implies that the string size $N$ of the most recent measurements and actions in $s^n$ and $s$ that are identical goes to $\infty$ as $n \to \infty$ and is upper bounded so that $\bar{d}(s^n,s) < \epsilon \implies N > \log_{\rho}(\frac{1}{\epsilon})$. Under mixing conditions, \cite[Assumption 5.10]{demirciRefined2023}, this then implies that, by \cite[Corollary 4.2]{le2004stability} (see also \cite[Lemma 5.11]{demirciRefined2023}) and the proof of \cite[Lemma 6]{demirci2023geometric}, $\lim_{N \to \infty} P(X_n \in \cdot | y_{(-N+1,0]},u_{(-N+1,-1]}, \pi_{-N} = \pi^*)$ is well-defined (as a limit of a Cauchy sequence), and this limit does not depend on $\pi^* \in {\cal P}(\mathbb{X})$ and only does so on $s$. This analysis implies then that $P(X_0 \in \cdot | y_{(-\infty,0]},u_{(-\infty,-1]})$ is well-defined and not only measurable but also continuous in $y_{(-\infty,0]},u_{(-\infty,-1]}$.
 
 {\it Wasserstein Regularity.} We finally present a sufficient condition for Assumption \ref{Wasser1}.  As in the above, we have that if $\bar{d}(s^n,s) < \epsilon$ then
\begin{align}
W_1(P(X_0 \in \cdot | (y,u)^m_{(-\infty,-N]},y_{(-N+1,0]},u_{(-N+1,-1]} ), P(X_0 \in \cdot | y_{(-\infty,0]},u_{(-\infty,-1]}) \nonumber \\
 \leq 2 \tau^{-\log_{\rho}(\frac{1}{\epsilon})} = 2 \tau^{\log_{\rho}(\epsilon)} , \nonumber
 \end{align}
where $\tau < 1$ is a contraction coefficient obtained via the Hilbert projective metric. Taking $\rho \in (\tau,1)$ ensures that 
$2 \tau^{\log_{\rho}(\epsilon)} \leq M_1 \epsilon$ for some constant $M_1 < \infty$, implying Assumption \ref{Wasser1}(i). Following as in \cite[Lemma 1]{demirciRefined2023}, and if we also impose Assumption \ref{main_assumption}-(\ref{CostLipschitz}) on the cost function, the Wasserstein regularity of the conditional probability as shown then implies Assumption \ref{Wasser1}(ii). 

%One can furthermore add regularity on the measurement
% Perhaps we can have $M_1 < 1$ also for the average cost problem.

\section{Proofs}

\subsection{Proof of Theorem \ref{Belief-MDPMethod}}\label{proofBelief-MDPMethod}

Following \cite{survey,Borkar2}, we study the limit distribution of the following occupation measures, under any policy $\gamma$ in $\Gamma$. Let for $T \geq 1$
\[v_T(D) = {1 \over T} \sum_{t=0}^{T-1} 1_{\{(\pi_t,u_t) \in D\}}, \quad D \in {\cal B}({\cal P}(\mathbb{X}) \times \mathbb{U}).\]
Consider any policy $\gamma$ in $\Gamma$, $\pi_0 \sim \mu$, and let for $T \geq 1$, \[\mu_T(D)  = E_{\mu}^{\gamma}[v_T(D)] =  E^{\gamma}_{\mu} {1 \over T} \bigg[ \sum_{t=0}^{T-1} 1_{\{\pi_t,u_t) \in D\}} \bigg], \quad D \in {\cal B}({\cal P}(\mathbb{X}) \times \mathbb{U})\]
Let for $\mu \in {\cal P}(\mathbb{X})$, $\mu P(A) := \int \mu(d\pi,du) P(\pi_{t+1} \in A | \pi_t=\pi, u_t=u)$. Then through what is often referred to as a {\it Krylov-Bogoliubov-type} argument, for every Borel $A \subset {\cal P}(\mathbb{X})$
\begin{eqnarray}
 |\mu_{T}(A \times \mathbb{U}) - \mu_{T} P (A)| = E^{\gamma}_{\mu_0} {1 \over T} \bigg[ \sum_{t=0}^{T-1} 1_{\{\pi_t,u_t) \in (A \times \mathbb{U})\}} - \sum_{t=0}^{T-1} 1_{\{\pi_{t+1},u_{t+1}) \in (A \times \mathbb{U})\}} \bigg] \leq {1 \over T} \to 0, \nonumber %\label{converEmp}
 \end{eqnarray}
 as $T \to \infty$. 
 Notice that the above applies for any policy $\gamma \in \Gamma$. Now, if we can ensure that, for some subsequence $\mu_{t_k}$, $\mu_{t_k} \to \mu$ weakly for some probability measure $\mu$ (which holds by the assumption that the state and action spaces are compact, which makes the set of probability measures on the state weakly compact), it would follow that $\mu_{t_k}\eta(A \times \mathbb{U}) \to \mu \eta(A \times \mathbb{U})$. It would follow that $\mu_{t_k}P \to \nu P$ also, by the weak Feller condition, and hence $\nu=\nu P $ and $\nu$ would be stationary.  

%Now define
%\begin{eqnarray}\label{ergodInv1}
%{\cal G}_{{\cal P}(\mathbb{X})} = \{v \in {\cal P}({\cal P}(\mathbb{X}) \times \mathbb{U}): &&v(B \times \mathbb{U}) = \int_{\pi,u} P(\pi_{t+1} \in B|\pi,u) v(d\pi,du),\nonumber \\
%&& \quad\quad \quad \quad \quad\quad B \in {\cal B}({\cal P}(\mathbb{X})) \} 
%\end{eqnarray}

Define, with $\Gamma_S$ denoting the set of stationary control policies mapping the belief state to actions,
\begin{eqnarray}\label{ergodInv2}
{\cal G} = \{v \in {\cal P}({\cal P}(\mathbb{X}) \times \mathbb{U}): && \exists \gamma \in \Gamma_S, v(A) = \int_{\pi,u} P^{\gamma}((\pi_{t+1},u_{t+1}) \in A|x) v(dx,du), \nonumber \\
&& \quad \quad\quad \quad \quad \quad\quad \quad A \in {\cal B}({\cal P}(\mathbb{X}) \times \mathbb{U}) \} 
\end{eqnarray}

%It is evident that ${\cal G} \subset {\cal G}_{\mathbb{X}}$ since there are fewer restrictions for ${\cal G}_{\mathbb{X}}$. However, these two sets are indeed equal. In particular, for $v \in {\cal G}_{{\cal P}(\mathbb{X})}$, if we write: $v(A,B)=\int_{A} v(d\pi) v(u \in B|\pi)$, then, we can construct a consistent $v \in {\cal G}$: $v(B \times C) = \int_{\pi \in B} v(d\pi) v(u \in C|\pi)$. 

Thus, every weakly converging sequence $\mu_t$ will satisfy the above equation, and therefore, under any admissible policy, every converging occupation measure sequence converges to the set ${\cal G}$. Let us define
\[\gamma^* = \inf_{v \in {\cal G}} \int v(d\pi,du) \tilde{c}(\pi,u)\]

Since the expression $\int v(dx,du) c(x,u)$ is lower semi-continuous in $v$, a compactness condition on ${\cal G}$ will ensure the existence of an optimal occupation measure which is in ${\cal G}$: we now show that there exists an optimal occupation measure in ${\cal G}$ if the transition kernel is weak-Feller and ${\cal G}$ is weakly compact: The problem has now reduced to
$$\inf_{\mu \in {\cal G}} \int \mu(dx,du) c(x,u), $$
The set ${\cal G}$ is closed, since if $\nu_n \to \nu$ and $\nu_n \in {\cal G}$, then for continuous and bounded $f \in C_b({\cal P}(\mathbb{X}))$, $\langle \nu_n, f \rangle \to \langle \nu, f \rangle$. By weak-Feller continuity of the kernel $\int f(\pi')\eta(d\pi'|\pi,u)$ is also continuous and thus, $\langle \nu_n, \eta f \rangle \to \langle \nu, \eta f \rangle = \langle \nu \eta, f \rangle$. Thus, $\nu(f) = \nu \eta(f)$ and $\nu \in {\cal G}$. Therefore, ${\cal G}$ is weakly sequentially compact. Since the integral $ \int \mu(dx,du) c(x,u)$ is lower semi-continuous on the set of measures under weak convergence, and the existence result follows from Weierstrass' Theorem. As a result, there exists an optimal occupation measure, say $v(d\pi, du)$. This defines a stationary control policy by the Radon-Nikodym derivative: $\mu(u \in \cdot | \pi) = {d v(d\pi, \cdot) \over d \int_{u} v(d\pi, \cdot) }(\pi)$, for $v$ a.e. $\pi$.
\qed

\subsection{Proof of Lemma \ref{keyLemmaClosed}}\label{Proofkeylemmaclosed}
%\textbf{Proof.}\{of Lemma \ref{keyLemmaClosed} \}
a) Let us recall the $w$-$s$ topology \cite{balder2001,Schal}: Let $\mathbb{A}, \mathbb{B}$ be complete, separable, metric spaces. The $w$-$s$ topology on the set of probability measures ${\cal P}(\mathbb{A} \times \mathbb{B})$ is the coarsest topology under which $\int f(a,b) \nu(da,db): {\cal P}(\mathbb{A} \times \mathbb{B}) \to \mathbb{R}$ is continuous for every measurable and bounded $f$ which is continuous in $b \in \mathbb{B}$ for every $a \in \mathbb{A}$ (but unlike weak topology, $f$ does not need to be continuous in $a$).

For every $n$, we have that for every function $g$ that is continuous and bounded:
\begin{eqnarray}\label{conditionalIndepEqn}
&& \int P_n(dy^1,dy^2,du^2)g(y^1,y^2,u^2) = \int P_n(dy^1|y^2) P_n(du^2,y^2)g(y^1,y^2,u^2) \nonumber \\
&& = \int P(dy^1|y^2) P_n(du^2,y^2)g(y^1,y^2,u^2) 
\end{eqnarray}

Testing the equality above on continuous and bounded functions implies this property for any measurable and bounded function (that is, continuous and bounded functions form a {\it separating class}, see e.g. p. 13 in \cite{Bil99} or Theorem 3.4.5 in \cite{ethier2009markov}) for weak convergence of probability measures. 
Since the marginals on $y^1,y^2$ is fixed, \cite[Theorem 3.10]{Schal} (see also \cite[Theorem 2.5]{balder2001}) establishes that the sequence $\{P_n\}$ is relatively compact under the $w$-$s$ topology under the stated tightness condition; $\{P_n\}$ is tight by Prohorov's theorem.

Now, taking the limit of both sides in (\ref{conditionalIndepEqn}), we have that the left hand side converges to $\int P(dy^1,dy^2,du^2)g(y^1,y^2,u^2)$. The right-hand side, on the other hand can be written as:
\begin{eqnarray}
\int \bigg(\int P(dy^1|dy^2) g(y^1,y^2,u^2) \bigg) P_n(du^2,dy^2) 
\end{eqnarray}
The expression $\bigg(\int P(dy^1|dy^2) g(y^1,y^2,u^2) \bigg)$ is measurable in $y^2$ and continuous in $u^2$ by an application of the dominated convergence theorem. As a result, by the $w$-$s$ convergence, in the limit we have
\begin{eqnarray}
\int \bigg(\int P(dy^1|dy^2) g(y^1,y^2,u^2) \bigg) P(du^2,y^2) 
\end{eqnarray}
and thus the conditional independence property is satisfied.

%\subsection{Proof of Lemma \ref{keyLemmaClosedCont}}\label{Proofkeylemmacontrolclosed}

%\textbf{Proof.}\{of Lemma \ref{keyLemmaClosedCont} \}
b) For every $n$, we again have that for every function $g$ continuous and bounded:
\begin{eqnarray}\label{conditionalIndepEqn1}
&& \int P_n(dy^1,dy^2,du^2)g(y^1,y^2,u^2) = \int P_n(dy^1|y^2) P_n(du^2,y^2)g(y^1,y^2,u^2) \nonumber \\
&& = \int \kappa(dy^1|y^2) P_n(du^2,y^2)g(y^1,y^2,u^2) 
\end{eqnarray}
The proof then closely follows the arguments above in a) and is omitted for space constraints.

%As in (i), testing the equality above on continuous and bounded functions implies this property for any measurable and bounded function \cite[p. 13]{Bil99} or \cite[Theorem 3.4.5]{ethier2009markov}) for weak convergence of probability measures. 
%Since the marginals on $y^1,y^2$ is fixed, \cite[Theorem 3.10]{Schal} (see also \cite[Theorem 2.5]{balder2001}) establishes that the sequence $\{P_n\}$ is relatively compact under the $w$-$s$ topology under the stated tightness condition; $\{P_n\}$ is tight by Prohorov's theorem. Now, taking the limit of both sides in (\ref{conditionalIndepEqn1}), we have that the left hand side converges to $\int P(dy^1,dy^2,du^2)g(y^1,y^2,u^2)$. The right-hand side can be written as:
%\begin{eqnarray}
%\int \bigg(\int \kappa(dy^1|dy^2) g(y^1,y^2,u^2) \bigg) P_n(du^2,y^2) 
%\end{eqnarray}
%The expression $\bigg(\int \kappa(dy^1|dy^2) g(y^1,y^2,u^2) \bigg)$ is measurable in $y^2$ and continuous in $u^2$ by an application of dominated convergence theorem. Furthermore, by the condition that $\kappa$ is a continuous kernel, we have that $\bigg(\int \kappa(dy^1|dy^2) g(y^1,y^2,u^2) \bigg)$ is continuous in $y^2$, by an application of a generalized dominated convergence theorem \cite[Theorem 3.5]{serfozo1982convergence} or \cite[Theorem 3.5]{Lan81}. As a result, in the limit we have
%\begin{eqnarray}
%\int \bigg(\int P(dy^1|dy^2) g(y^1,y^2,u^2) \bigg) P(du^2,y^2) 
%\end{eqnarray}
%and thus the conditional independence property is satisfied.
\qed

%\qed

%\sy{We can update our SICON paper. No need for weak continuity for the existence result for classical information structures.!!!!!! Update this as such at least for the book's second edition. WE SHOULD.}
%
%\sy{However, in the above it is crucial that we have the stationary measure for the measurement processes!}

\subsection{Proof of Theorem \ref{controlledPOMDPE}}\label{controlledPOMDPEProofSection}

Recall (\ref{ergodInv2}) and let \[\mu_T(D)  = E[v_T(D)] =  E_{v_0} {1 \over T} \bigg[ \sum_{t=1}^T 1_{\{z_t,u_t) \in D\}} \bigg], \quad D \in {\cal B}(\mathbb{X} \times \mathbb{Y}^{\mathbb{Z}_-} \times \mathbb{U}^{\mathbb{Z}_-} \times \mathbb{U}).\]

This can also be written as ${1 \over T} \bigg[ \sum_{t=1}^T P( \{(z_t,u_t) \in D\}) \bigg]$. Recall that (\ref{CSIP}) holds. Let $(x_t,y_{(-\infty,t]},u_t) \sim P_t$.  As in the proof of Theorem \ref{proofBelief-MDPMethod}, through a {\it Krylov-Bogoliubov-type} argument, for every Borel $A \in {\cal B}(\mathbb{X} \times \mathbb{Y}^{\mathbb{Z}_-} \times \mathbb{U}^{\mathbb{Z}_-}$
\begin{eqnarray}
&& |\mu_{N}(A \times \mathbb{U}) - \mu_{N} \mathbb{P}^a (A)| \nonumber \\
&& = \bigg| {1 \over N} \bigg( (\mu_0(A \times \mathbb{U})+\cdots+ \mu (\mathbb{P}^a)^{(N-1)}(A))  - (mu_0 \mathbb{P}^a(A)+\cdots+v (\mathbb{P}^a)^{N}(A))  \bigg) \bigg| \nonumber\\
&& \leq {1 \over N} |\mu_0(A \times \mathbb{U}) - \mu_0 (\mathbb{P}^a)^{N}(A)|
 \to 0. \label{converEmp}
 \end{eqnarray}
As earlier, in the proof of Theorem \ref{proofBelief-MDPMethod}, if we can ensure that for some subsequence, $\mu_{t_k} \to \mu$ for some probability measure $\mu$, it would follow that $\mu_{t_k}P \to \mu$ also and hence $\mu=\mu P$ and $\mu$ would be invariant and $\mu \in {\cal G}$, where
\begin{eqnarray}\label{ergodInv33}
{\cal G} &=& \{v \in {\cal P}(\mathbb{X} \times \mathbb{Y}^{\mathbb{Z}_-} \times \mathbb{U}^{\mathbb{Z}_-}): \nonumber \\
&& v(B \times \mathbb{U}) = \int_{x,u} \mathbb{P}^a(z_{t+1} \in B|z,u) v(dz,du), B \in {\cal B}(\mathbb{X} \times \mathbb{Y}^{\mathbb{Z}_-} \times \mathbb{U}^{\mathbb{Z}_-}) \}  \nonumber 
\end{eqnarray}
Convergence to ${\cal G}$ is by (\ref{converEmp}) is a consequence of Lemma \ref{weakFellerP} under Assumption \ref{contCond}. 

We will show that, if convergence occurs, it must also be that $\mu \in {\cal H}$;
\begin{eqnarray}\label{ergodInv44}
{\cal H} &=& \{v \in {\cal P}(\mathbb{X} \times \mathbb{Y}^{\mathbb{Z}_-} \times \mathbb{U}^{\mathbb{Z}_-}): \nonumber \\
&& \quad  v(dx_0,dy_{(-\infty,0]},du_{(-\infty,-1]},du_0) \nonumber \\
&&  \quad \quad = v(dx_0,dy_{(-\infty,0]},du_{(-\infty,-1]}) v(du_0 | dy_{(-\infty,0]},du_{(-\infty,-1]})\}, 
\end{eqnarray}
that is the control action variables are conditionally independent from the state variables given the information variable $y_{(-\infty,0]},u_{(-\infty,-1]}$, thus satisfying (\ref{CSIP}). 

We show now that $\mu_t$, possibly along a subsequence, will also converge to ${\cal H}$. Observe that
\begin{eqnarray}\label{condIndP222}
&& \int \mu_T(dx,dy_{(-\infty,0]},du_{(-\infty,-1]},du) f(x,y_{(-\infty,0]},u_{(-\infty,-1]},u) \nonumber \\
&& =  \int {1 \over T}  \sum_{t=1}^T P_t(dx,dy_{(-\infty,0]},du_{(-\infty,-1]},du_0) f(x,y_{(-\infty,0]},u_{(-\infty,-1]},u_0) \nonumber \\
&& =  \int {1 \over T}  \sum_{t=1}^T P_t(dx| y_{(-\infty,0]},du_{(-\infty,-1]}) P_t(dy_{(-\infty,0]},du_{(-\infty,-1]},du_0) f(x,y_{(-\infty,0]},u_{(-\infty,-1]},u_0) \nonumber \\
&& =  \int {1 \over T}  \sum_{t=1}^T \bigg( P_t(dx| y_{(-\infty,0]},du_{(-\infty,-1]}) f(x,y_{(-\infty,0]},u_{(-\infty,-1]},u_0)  \bigg) P_t(dy_{(-\infty,0]},du_{(-\infty,-1]},du_0) \nonumber \\
&& =  \int {1 \over T}  \sum_{t=1}^T \bigg( P(dx| y_{(-\infty,0]},du_{(-\infty,-1]}) f(x,y_{(-\infty,0]},u_{(-\infty,-1]},u_0)  \bigg) P_t(dy_{(-\infty,0]},du_{(-\infty,-1]},du_0) \nonumber \\
\label{stat1112} \\
&& =  \int {1 \over T}  \sum_{t=1}^T \bigg( g(y_{(-\infty,0]},u_{(-\infty,-1]},u_0) \bigg) P_t(dy_{(-\infty,t]},du_{(-\infty,t-1]},du_t) \nonumber \\ \label{stat22} \\
&& =  \int {1 \over T}  \sum_{t=1}^T \bigg( g(y_{(-\infty,t]},u_{(-\infty,t-1]},u_t) \bigg) P_t(dy_{(-\infty,t]},du_{(-\infty,t-1]},du_t) \nonumber \\
&& \to \int \bigg( g(y_{(-\infty,0]},u_{(-\infty,-1]},u_0) \bigg) \mu(dy_{(-\infty,0]},du_{(-\infty,-1]},du_0)  \label{weakcontG2} \\
&& =  \int \bigg( P(dx| y_{(-\infty,0]},u_{(-\infty,-1]}) f(x,y_{(-\infty,0]},u_{(-\infty,-1]},u_0)  \bigg)  \mu(dy_{(-\infty,0]},du_{(-\infty,-1]},du_0)  \nonumber \\
&& =  \int \bigg( P(dx| y_{(-\infty,0]},u_{(-\infty,-1]}) \mu(dy_{(-\infty,0]},du_{(-\infty,-1]},du_0) \bigg)    f(x,y_{(-\infty,0]},u_{(-\infty,-1]},u_0)  
\end{eqnarray}
Here (\ref{stat1112}) is due to Assumption \ref{ForgetCondition}\footnote{In the absence of Assumption \ref{ForgetCondition} (\ref{stat1112}) would be incorrect; as a counterexample, consider $y_k$ giving no information (e.g. $y_k$ being constant); in this case, we clearly require $X_t$ to be stationary for this to be correct since there is no information: $P(X_t | Y_{-\infty,t]}) = P(X_{t+1} | Y_{-\infty,t+1]})$.} and in (\ref{stat22}) we define \[g(y_{(-\infty,0]},u_{(-\infty,-1]},u_0) = \bigg( P(dx| y_{(-\infty,0]},u_{(-\infty,-1]}) f(x,y_{(-\infty,0]},u_{(-\infty,-1]},u_0)  \bigg).\]
Here, (\ref{weakcontG2}) follows from weak continuity of the kernel under Assumption \ref{contConditinal}. Thus, if $\mu_t(\cdot)$ converges weakly to $\mu(\cdot)$, $\mu$ satisfies the conditional independence property and is thus in ${\cal H}$. 

In view of this, the existence problem reduces to 
\begin{eqnarray}\label{expCostGH}
\inf_{v \in {\cal G} \cap {\cal H}} \int v(dx,dy_{(-\infty,0]},du_{(-\infty,-1]},du) c(x,u)
\end{eqnarray}

By Lemma \ref{closednessTight2} we have that ${\cal G} \cap {\cal H}$ is a compact set under weak topology, and since $c$ is continuous, there exists an optimal measure $v \in {\cal G} \cap {\cal H}$.

\qed

\subsection{Proof of Theorem \ref{FilterInvarianceMeasure2}}

Let $\bar{P}$ instead denote an invariant process measure on $(Y_{(-\infty,k]},U_{(-\infty,k-1]}, U_k, X_k)$. Let $P_0$ be the projection of $\bar{P}$ on $dy_{(-\infty,-1]},du_{(-\infty,-1]}$. If $P_0(dy_{(-\infty,-1]},du_{(-\infty,-1]}) \times \pi_0(dx_0,dy_0,du_u) \ll \bar{P}$; that is, with the initial state distribution selected independently of the {\it past} process; then with $u_t=\gamma(y_{(-\infty,t]},u_{(-\infty,t-1]},r_t)$ and \[g(y_{(-\infty,k]},u_{(-\infty,k-1]},x_k) = E[c(x_k,u_k) |x_k, y _{(-\infty,k]},u_{(-\infty,k-1]}],\]
\[ \lim_{T \to \infty} {1 \over T} E_{P_0 \times \pi_0} \sum_{k=0}^{T-1} g(y_{(-\infty,k]},u_{(-\infty,k-1]},x_k) = \int g(y_{(-\infty,k]},u_{(-\infty,k-1]},x_k) \bar{Q}(y_{(-\infty,0]},u_{(-\infty,-1]}) \]
for some $\bar{Q}$ which is invariant. Note here that $\bar{Q}$ does not need to be equal to $\bar{P}$.

%Let $u_k = \gamma(y_{(-\infty,k]})$ be optimal (one may also add independent randomization). In this case, define $g(y_{(-\infty,k]},x_k) = c(x_k, u_k)$. Thus,
%$P_0 \times \pi_0$ a.s.
%\[ \lim_{T \to \infty} {1 \over T} E_{y_{(-\infty,0},x_0} [\sum_{k=0}^{T-1} g(y_{(-\infty,k]},x_k)] = \int g(y_{(-\infty,k]},x_k) \bar{Q}(dy_{(-\infty,0]}),\]
%by the absolute continuity condition since the set of initial conditions which does not satisfy convergence has zero measure under $\bar{P}$ and thus under $P_0 \times \pi_0$.
%
%But we know that an optimal control policy will use the conditional probability measure to select the action at any given time. Thus, we can write $u_k = \gamma(y_{(-\infty,k]}) = \gamma(\pi_k)$ with $\pi_k(\cdot) = P(x_k \in \cdot | y_{(-\infty,k})$. Since $P_0$ is stationary, we have that there exists a stationary measure on $\pi_t$ given by:
%\[\eta(\pi \in B) = P_0(F(Y_{(-\infty,-1]}) \in B) ),\]
%so that $\eta(B)= P_0(F^{-1}(B))$ for any Borel $B \in {\cal P}(\mathbb{X})$ under the weak convergence topology on ${\cal P}(\mathbb{X})$.
%Thus, instead of generating $y_{(-\infty,-1]}$ randomly according to $P_0$, we can generate $\pi_{-1}$ according to $\eta$. Convergence will hold with probability $1$.  
\qed

\subsection{Proof of Theorem \ref{MeanOptimal}}\label{SMeanOptimalProof}

The proof closely follows that of Theorem \ref{controlledPOMDPE}, and builds on the convergence properties of the sequence $\mu_T(D)  = E[v_T(D)] =  E_{v_0} {1 \over T} \bigg[ \sum_{t=1}^T 1_{\{(z_t,u_t) \in D\}} \bigg], \quad D \in {\cal B}(\mathbb{X} \times \mathbb{Y}^{\mathbb{Z}_-} \times \mathbb{U}^{\mathbb{Z}_-} \times \mathbb{U})$. We only note that, convergence to ${\cal H}$, under stationarity (Assumption \ref{statX}) is more direct. %The detailed proof is omitted due to space constraints and is available in \cite{anotherLookPOMDPs}. 

Recall (\ref{ergodInv2}) and let \[\mu_T(D)  = E[v_T(D)] =  E_{v_0} {1 \over T} \bigg[ \sum_{t=1}^T 1_{\{(z_t,u_t) \in D\}} \bigg], \quad D \in {\cal B}(\mathbb{X} \times \mathbb{Y}^{\mathbb{Z}_-} \times \mathbb{U}^{\mathbb{Z}_-} \times \mathbb{U}).\]

This can also be written as
\[ {1 \over T} \bigg[ \sum_{t=1}^T P( \{(z_t,u_t) \in D\}) \bigg] \]
Observe that for any $t$, $x_t \leftrightarrow y_{(-\infty,t]} \leftrightarrow u_t$ holds. Let $(x_t,y_{(-\infty,t]},u_t) \sim P_t$. Now, again through the {\it Krylov-Bogoliubov-type} argument, for every Borel $A$
\begin{eqnarray}
&& |\mu_{N}(A \times \mathbb{U}) - \mu_{N} \mathbb{P}^a (A)| \nonumber \\
&& = \bigg| {1 \over N} \bigg( (\mu_0(A  \times \mathbb{U})+\cdots+ \mu (\mathbb{P}^a)^{(N-1)}(A))  - (v_0 \mathbb{P}^a(A)+\cdots+v (\mathbb{P}^a)^{N}(A))  \bigg) \bigg| \nonumber\\
&& \leq {1 \over N} |\mu_0(A  \times \mathbb{U}) - \mu_0 (\mathbb{P}^a)^{N}(A)|
 \to 0. \label{converEmp}
 \end{eqnarray}
Now, if we can ensure that for some subsequence, $\mu_{t_k} \to \mu$ for some probability measure $\mu$, it would follow that $\mu_{t_k}P(B) \to \mu(B)$ also and hence $\mu(B)=\mu P (B)$ and $\mu$ would be invariant.  

Thus, if $\mu_t$ converges, it must be that $\mu_t \to \mu$ for some $\mu \in  {\cal G} \cap {\cal H}$, where
\begin{eqnarray}\label{ergodInv3}
{\cal G} &=& \{v \in {\cal P}(\mathbb{X} \times \mathbb{Y}^{\mathbb{Z}_-} \times \mathbb{U}):  \nonumber \\
&& v(B \times \mathbb{U}) = \int_{x,u} P(X_1,Y_{(-\infty,1]} \in B|x_0,y_{(-\infty,0]},u) v(dx_0,dy_{(-\infty,0]},du), \nonumber \\
&& \quad \quad  \quad \quad  \quad \quad  \quad \quad  \quad \quad \quad B \in {\cal B}(\mathbb{X} \times \mathbb{Y}^{\mathbb{Z}_-} \times \mathbb{U}) \} 
\end{eqnarray}
and
\begin{eqnarray}\label{ergodInv4}
{\cal H} &=& \{v \in {\cal P}(\mathbb{X} \times \mathbb{Y}^{\mathbb{Z}_-} \times \mathbb{U}): \nonumber \\
&&  v(dx,dy_{(-\infty,0]},du_0) = v(dx,dy_{(-\infty,0]}) P(du_0 | dy_{(-\infty,0]})\}, 
\end{eqnarray}
that is the control action variables are conditionally independent from the state variables given the information variable $y_{(-\infty,0]},u_{(-\infty,-1]}$.

Convergence to {\cal G} is by (\ref{converEmp}) as a consequence of Lemma \ref{weakFellerP} under Assumption \ref{contCond}. 

We show now that $\mu_t$ will also converge to ${\cal H}$. 
% is by (\ref{condIndP}) and (\ref{condIndP2}) below.

{\bf Case 1. Under Assumption \ref{statX}.}
For a continuous and bounded $f$, if $\mu_t \to \mu$ for some $\mu$, we have that
\begin{eqnarray}\label{condIndP}
&& \int \mu_t(dx,dy_{(-\infty,0]},du) f(x,y_{(-\infty,0]},u) \nonumber \\
&& =  \int {1 \over T}  \sum_{t=1}^T P_t(dx,dy_{(-\infty,0]},du_0) f(x,y_{(-\infty,0]},u_0) \nonumber \\
&& =  \int {1 \over T}  \sum_{t=1}^T P_t(dx| y_{(-\infty,0]}) P_t(dy_{(-\infty,0]},du_0) f(x,y_{(-\infty,0]},u_0) \nonumber \\
&& =  \int {1 \over T}  \sum_{t=1}^T \bigg( P_t(dx| y_{(-\infty,0]}) f(x,y_{(-\infty,0]},u_0)  \bigg) P_t(dy_{(-\infty,0]},du_0) \nonumber \\
&& =  \int {1 \over T}  \sum_{t=1}^T \bigg( P(dx| y_{(-\infty,0]}) f(x,y_{(-\infty,0]},u_0)  \bigg) P_t(dy_{(-\infty,0]},du_0) \label{stat11} \\
&& =  \int {1 \over T}  \sum_{t=1}^T \bigg( g(y_{(-\infty,0]},u_0) \bigg) P_t(dy_{(-\infty,t]},du_t) \label{stat2} \\
&& =  \int {1 \over T}  \sum_{t=1}^T \bigg( g(y_{(-\infty,t]},u_t) \bigg) P_t(dy_{(-\infty,t]},du_t) \nonumber \\
&& \to \int \bigg( g(y_{(-\infty,0]},u_0)) \bigg) \mu(dy_{(-\infty,0]},du_0) \nonumber \\
&& =  \int \bigg( P(dx| y_{(-\infty,0]}) f(x,y_{(-\infty,0]},u_0)  \bigg)  \mu(dy_{(-\infty,0]},du_0)  \nonumber \\
&& =  \int \bigg( P(dx| y_{(-\infty,0]}) \mu(dy_{(-\infty,0]},du_0)\bigg)    f(x,y_{(-\infty,0]},u_0)  
\end{eqnarray}
where (\ref{stat11}) uses stationarity and in (\ref{stat2}) we define \[g(y_{(-\infty,0]},u_0) = \bigg( P(dx| y_{(-\infty,0]}) f(x,y_{(-\infty,0]},u_0)  \bigg).\]
This expression will converge as soon as $\mu_t(dy_{(-\infty,0]},du_0)$ converges weakly to $\mu$, where $\mu$ satisfies the conditional independence property. Here, $\mu_t$ converges to $\mu$ in the $w-s$ sense. Thus, even in the absence of Assumption \ref{contConditinal}, convergence holds in this case. \footnote{Here, we crucially assume that the marginal on $x_k, y_{(-\infty,k]}$ is fixed, and the realizations are generated according to the stationary measure. This is crucial for the argument in (\ref{condIndP}). In the above, we don't need weak continuity.}

%\footnote{A cautious reader will argue that the optimal policy is to pick $\arg \min E[c(X_k,U_k)|Y_{(-\infty,k]}]$. This clearly is a stationary policy if the filter is %stationary.}

{\bf Case 2. Under Assumptions \ref{contConditinal} and \ref{ForgetCondition}.}
If we don't assume stationarity, we can modify the above as follows, but then we need the weak continuity condition stated in Assumption \ref{contConditinal}:

\begin{eqnarray}\label{condIndP22}
&& \int \mu_t(dx,dy_{(-\infty,0]},du) f(x,y_{(-\infty,0]},u) \nonumber \\
&& =  \int {1 \over T}  \sum_{t=1}^T P_t(dx,dy_{(-\infty,0]},du_0) f(x,y_{(-\infty,0]},u_0) \nonumber \\
&& =  \int {1 \over T}  \sum_{t=1}^T P_t(dx| y_{(-\infty,0]}) P_t(dy_{(-\infty,0]},du_0) f(x,y_{(-\infty,0]},u_0) \nonumber \\
&& =  \int {1 \over T}  \sum_{t=1}^T \bigg( P_t(dx| y_{(-\infty,0]}) f(x,y_{(-\infty,0]},u_0)  \bigg) P_t(dy_{(-\infty,0]},du_0) \nonumber \\
&& =  \int {1 \over T}  \sum_{t=1}^T \bigg( P(dx| y_{(-\infty,0]}) f(x,y_{(-\infty,0]},u_0)  \bigg) P_t(dy_{(-\infty,0]},du_0) \label{stat111} \\
&& =  \int {1 \over T}  \sum_{t=1}^T \bigg( g(y_{(-\infty,0]},u_0) \bigg) P_t(dy_{(-\infty,t]},du_t) \label{stat222} \\
&& =  \int {1 \over T}  \sum_{t=1}^T \bigg( g(y_{(-\infty,t]},u_t) \bigg) P_t(dy_{(-\infty,t]},du_t) \nonumber \\
&& \to \int \bigg( g(y_{(-\infty,0]},u_0)) \bigg) \mu(dy_{(-\infty,0]},du_0)  \label{weakcontG} \\
&& =  \int \bigg( P(dx| y_{(-\infty,0]}) f(x,y_{(-\infty,0]},u_0)  \bigg)  \mu(dy_{(-\infty,0]},du_0)  \nonumber \\
&& =  \int \bigg( P(dx| y_{(-\infty,0]}) \mu(dy_{(-\infty,0]},du_0)\bigg)    f(x,y_{(-\infty,0]},u_0)  
\end{eqnarray}
\footnote{In the absence of Assumption \ref{ForgetCondition} (\ref{stat111}) would be incorrect. As a counterexample, consider $y_k=$ no information; in this case, we clearly require $X_t$ to be stationary for this to be correct since there is no information: $P(X_t | Y_{-\infty,t]}) = P(X_{t+1} | Y_{-\infty,t+1]})$.}
Here (\ref{stat111}) is due to Assumption \ref{ForgetCondition} and in (\ref{stat222}) we define \[g(y_{(-\infty,0]},u_0) = \bigg( P(dx| y_{(-\infty,0]}) f(x,y_{(-\infty,0]},u_0)  \bigg).\]
This expression will converge as soon as $\mu_t(dy_{(-\infty,0]},du_0)$ converges weakly to $\mu$, where $\mu$ satisfies the conditional independence property. Here, (\ref{weakcontG}) follows from weak continuity of the kernel under Assumption \ref{contConditinal}.

In view of this, the existence problem reduces to 
\[\inf_{v \in {\cal G} \cap {\cal H}} \int v(dx,dy_{(-\infty,0]},du) c(x,u)\]
where

We note that under weak continuity of the transition kernel for $\mathbb{P}^a$, ${\cal G}$ is a closed set under weak convergence. By Lemma \ref{keyLemmaClosed}, ${\cal H}$ is also closed. Since the state space is compact, the set of probability measures on $(\mathbb{X} \times \mathbb{Y}^{\mathbb{Z}_-} \times \mathbb{U})$ is tight. Thus ${\cal G} \cap {\cal H}$ is a compact set under weak topology and since $c$ is continuous, there exists an optimal measure $v \in {\cal G} \cap {\cal H}$.
\qed

\section{An ergodic theorem for Markov chains}\label{appendixA}

Suppose that $\{X_t\}_{t\ge 0}$ denote a discrete-time Markov chain with state space
$\mathbb{X}$, a Polish space.
 \begin{theorem}\label{convergenceT3} \cite{HernandezLermaLasserre} \cite{worm2010ergodic}
Let $\bar{P}$ be an invariant probability measure for a Markov process.
\begin{itemize}
\item[(i)] [Ergodic decomposition and weak convergence] For $x$, $\bar{P}$ a.s., ${1 \over N} E_x[\sum_{t=0}^{N-1} 1_{\{x_n \in \cdot\}}] \to P_x( \cdot)$ weakly and $\bar{P}$ is invariant for $P_x(\cdot)$ in the sense that
     \[\bar{P}(B) = \int P_x(B) \bar{P}(dx)   \]
\item[(ii)] [Convergence in total variation] For all $\mu \in {\cal P}(\mathbb{X})$ which satisfies that $\mu \ll \bar{P}$ (that is, $\mu$ is
 absolutely continuous with respect to $\bar{P}$), there exists $v^*$ such that
\[ \| E_{\mu}[{1 \over N}\sum_{t=0}^{N-1} 1_{\{X_t \in \cdot\}}]  - v^*(\cdot) \|_{TV} \to 0. \]
\end{itemize}
\end{theorem}

\end{document}